\documentclass[11pt]{article}
\usepackage{amsmath}
\usepackage{amsthm}
\usepackage{amssymb}
\usepackage{graphicx}
\usepackage{multirow} 
\usepackage{tabularx}
\newcommand{\Z}{\mathbb{Z}}
\newcolumntype{R}{>{\raggedleft\arraybackslash}X}
\newcolumntype{L}{>{\raggedright\arraybackslash}X}
\usepackage[a4paper,margin=3cm,right=3cm,includehead]{geometry}
\theoremstyle{plain}
\newtheorem*{theorem}{Theorem 3}
\title{ \bf Full proof of the existence of a degree 8 circulant graph of order $L(8,k)$ of arbitrary diameter $k$}
\author{R. R. Lewis\\[-3pt]
\small Department of Mathematics and Statistics\\[-3pt]
\small The Open University\\[-3pt]
\small Milton Keynes, UK\\[-3pt]
\small \texttt{robert.lewis@open.ac.uk}}
\date{10th April 2014}



\makeatletter
\renewcommand\section{\@startsection {section}{1}{\z@}%
                                   {-2.5ex \@plus -1ex \@minus -.2ex}%
                                   {1.3ex \@plus.2ex}%
                                   {\normalfont\bf}}
\makeatother

\begin{document}
\maketitle

This is the full proof of Theorem 3 in the paper \textquotedblleft The degree-diameter problem for circulant graphs of degree 8 and 9" by the author \cite{Lewis}.
To avoid the paper being unduly long it includes only the exceptions for the orthant of $\textbf{v}_1$ for diameter $k \equiv 0 \pmod 2$ and for $k \equiv 1 \pmod 2$. In the version below the exceptions for all eight orthants for diameter $k \equiv 0$ and $k \equiv 1$ (mod 2) are included in full. This proof closely follows the approach taken by Dougherty and Faber in their proof of the existence of the degree 6 graph of order $DF(6,k)$ for all diameters $k \geq 2$ \cite{Dougherty}.

\begin{theorem}
For all $k\geq 2$, there is an undirected Cayley graph on four generators of a cyclic group which has diameter k and order $L(8,k)$, where
\[ L(8,k)=
\begin{cases}
(k^4+2k^3+6k^2+4k)/2 &\mbox{ if } k\equiv 0  \pmod 2\\
(k^4+2k^3+6k^2+6k+1)/2 &\mbox{ if } k\equiv 1  \pmod 2 
\end{cases}
\]
Moreover for $k \equiv 0 \pmod 2$ a generator set is
$\{1, (k^3+2k^2+6k+2)/2, (k^4+4k^2-8k)/4, (k^4+4k^2-4k)/4\}$, \\
and for $k\equiv 1 \pmod 2$,
$\{1, (k^3+k^2+5k+3)/2, (k^4+2k^2-8k-11)/4, (k^4+2k^2-4k-7)/4\}$.
\label{theorem:A}
\end{theorem}
\begin{proof}
We will show the existence of four-dimensional lattices $L_k \subseteq \Z^4$ such that $\Z^4/L_k$ is cyclic, $S_k+L_k=\Z^4$, where $S_k$ is the set of points in $\Z^4$ at a distance of at most $k$ from the origin under the $l^1$ (Manhattan) metric, and $\vert \Z^4 : L_k\vert = L(8,k)$ as specified in the theorem. Then, by Theorem 1 of \cite {Lewis}, the resultant Cayley graph has diameter at most $k$. 

Let $a=\begin {cases} k/2 & \mbox { for } k \equiv 0 \pmod 2\\ (k+1)/2 & \mbox { for } k \equiv 1 \pmod 2. \end {cases}$

For $k \equiv 0 \pmod 2$ let $L_k$ be defined by four generating vectors as follows:
\[
\begin{array}{rcl}
\textbf{v}_1&=&(-a-1,a+1,a,-a+1)\\
\textbf{v}_2&=&(a-1,a+1,a+1,-a)\\
\textbf{v}_3&=&(-a-1,-a+1,a+1,-a)\\
\textbf{v}_4&=&(-a,-a,a,a+1)
\end{array}
\] 

Then the following vectors are in $L_k$:
\[
\begin{array}{l}
-(2a^2+2a+1)\textbf{v}_1 + (2a^2+a+2)\textbf{v}_2 - (a+2)\textbf{v}_3 +\textbf{v}_4 = (4a^3+4a^2+6a+1, -1, 0, 0),\\
-(2a^3-1)\textbf{v}_1 + (2a^3-a^2+2a-2)\textbf{v}_2 - (a^2+a-1)\textbf{v}_3 +(a-1)\textbf{v}_4 = (4a^4+4a^2-4a, 0, -1, 0),\\
-2a^3\textbf{v}_1 + (2a^3-a^2+2a-1)\textbf{v}_2 - (a^2+a-1)\textbf{v}_3 +(a-1)\textbf{v}_4 = (4a^4+4a^2-2a, 0, 0, -1)\\
\end{array}
\]
Hence we have $\textbf{e}_2 =  (4a^3+4a^2+6a+1)\textbf{e}_1,  \textbf{e}_3 =  (4a^4+4a^2-4a)\textbf{e}_1$  and $\textbf{e}_4 =  (4a^4+4a^2-2a)\textbf{e}_1$ in $\Z^4/L_k$, and so $\textbf{e}_1$ generates $\Z^4/L_k$.

Also det $\left ( \begin{array} {c} \textbf{v}_1 \\ \textbf{v}_2 \\ \textbf{v}_3 \\ \textbf{v}_4 \end {array} \right ) = $ det
$\left (
\begin{array} {l r r r}
8a^4+8a^3+12a^2+4a & 0 & 0 & 0 \\
4a^3+4a^2+6a+1 & -1 & 0 & 0 \\
4a^4+4a^2-4a & 0 & -1 & 0 \\
4a^4+4a^2-2a & 0 & 0 & -1 
\end {array} \right ) 
\newline
\newline
= -(8a^4+8a^3+12a^2+4a) = -(k^4+2k^3+6k^2+4k)/2= -L(8,k)$, as in the statement of the theorem. 

Thus $\Z^4/L_k$ is isomorphic to $\Z_{L(8,k)}$ via an isomorphism taking $\textbf{e}_1, \textbf{e}_2, \textbf{e}_3, \textbf{e}_4$ to $1$, $4a^3+4a^2+6a+1, 4a^4+4a^2-4a,  4a^4+4a^2-2a$. As $a=k/2$ this gives the first generator set specified in the theorem: $\{1, (k^3+2k^2+6k+2)/2, (k^4+4k^2-8k)/4, (k^4+4k^2-4k)/4\}$.

Similarly for $k \equiv 1 \pmod 2$ let $L_k$ be defined by four generating vectors as follows:
\[
\begin{array}{rcl}
\textbf{v}_1&=&(-a+1,a+1,-a+1,a)\\
\textbf{v}_2&=&(a+1,a+1,-a+2,a-1)\\
\textbf{v}_3&=&(-a-1,a-1,a-1,-a)\\
\textbf{v}_4&=&(-a,a,a,a-1)
\end{array}
\]

In this case the following vectors are in $L_k$:
\newline
$-(2a^2+a+2)\textbf{v}_1 + (2a^2+2a+1)\textbf{v}_2 - a\textbf{v}_3 -\textbf{v}_4 = (4a^3-4a^2+6a-1, -1, 0, 0)$,
\newline
$-(2a^3-a^2-2a-2)\textbf{v}_1 + (2a^3-4a-1)\textbf{v}_2 - (a^2-a-1)\textbf{v}_3 -(a-1)\textbf{v}_4 = (4a^4-8a^3+8a^2-8a, 0, -1, 0)$,
\newline
$-(2a^3-a^2-2a-1)\textbf{v}_1 + (2a^3-4a)\textbf{v}_2 - (a^2-a-1)\textbf{v}_3 -(a-1)\textbf{v}_4 = (4a^4-8a^3+8a^2-6a, 0, 0, -1)$.

Hence we have $\textbf{e}_2 = (4a^3+4a^2+6a-1)\textbf{e}_1$,  $\textbf{e}_3 = (4a^4-8a^3+8a^2-8a)\textbf{e}_1$ and $\textbf{e}_4=  (4a^4-8a^3+8a^2-6a)\textbf{e}_1$, in $\Z^4/L_k$, and so $\textbf{e}_1$ generates $\Z^4/L_k$.

Also det $\left ( \begin{array} {c} \textbf{v}_1 \\ \textbf{v}_2 \\ \textbf{v}_3 \\ \textbf{v}_4 \end {array} \right ) = $ det
$\left ( \begin{array} {l r r r}
8a^4-8a^3+12a^2-4a & 0 & 0 & 0 \\
4a^3-4a^2+6a-1 & -1 & 0 & 0 \\
4a^4-8a^3+8a^2-8a & 0 & -1 & 0 \\
4a^4-8a^3+8a^2-6a & 0 & 0 & -1 \end {array} \right ) 
\newline
\newline
= -(8a^4-8a^3+12a^2-4a) = -(k^4+2k^3+6k^2+6k+1)/2= -L(8,k)$, as in the statement of the theorem.

Thus $\Z^4/L_k$ is isomorphic to $\Z_{L(8,k)}$ with generators $1, 4a^3-4a^2+6a-1, 4a^4-8a^3+8a^2-8a, 4a^4-8a^3+8a^2-6a$.
As $a=(k+1)/2$ in this case, this gives the second generator set specified in the theorem: $\{1, (k^3+k^2+5k+3)/2, (k^4+2k^2-8k-11)/4, (k^4+2k^2-4k-7)/4\}$.

It remains to show that $S_k+L_k=\Z^4$. First we consider the case $k\equiv 0 \pmod 2$.

For $k=2$, it is straightforward to show directly that $\Z_{32}$ with generators $1, 4, 6, 15$ has diameter 2. So we assume $k\geq4$, so that $a \geq 2$. Now let
\[
\begin{array}{rcl}
\textbf{v}_5&=\textbf{v}_1-\textbf{v}_3+\textbf{v}_4=&(-a,a,a-1,a+2)\\
\textbf{v}_6&=\textbf{v}_1-\textbf{v}_2-\textbf{v}_4=&(-a,a,-a-1,-a)\\
\textbf{v}_7&=\textbf{v}_1-\textbf{v}_2-\textbf{v}_3=&(-a+1,a-1,-a-2,a+1)\\
\textbf{v}_8&=\textbf{v}_2-\textbf{v}_3+\textbf{v}_4=&(a,a,a,a+1)
\end{array}
\]
with $\textbf{v}_1, \textbf{v}_2, \textbf{v}_3, \textbf{v}_4$ as defined for $k \equiv 0$ (mod 2). Then the 16 vectors $\pm\textbf{v}_i$ for $i=1,...,8$ provide one element of $L_k$ lying strictly within each of the 16 orthants of $\Z^4$. Most of the coordinates of these vectors have absolute value at most $a+1$. Only $\pm \textbf{v}_5$ and $\pm \textbf{v}_7$ each have one coordinate with absolute value equal to $a+2$.

Now we consider the case $k\equiv 1 \pmod 2$. For $k=3$ it may be shown directly that $\Z_{104} $ with generators $1, 16, 20, 27$ has diameter 3. So we assume $k \geq5$, so that $a \geq3$, and let 
\[
\begin{array}{rcl}
\textbf{v}_5&=\textbf{v}_1-\textbf{v}_2-\textbf{v}_4=&(-a,-a,-a-1,-a+2)\\
\textbf{v}_6&=\textbf{v}_2+\textbf{v}_3-\textbf{v}_4=&(a,a,-a+1,-a)\\
\textbf{v}_7&=\textbf{v}_1+\textbf{v}_3-\textbf{v}_4=&(-a,a,-a,-a+1)\\
\textbf{v}_8&=\textbf{v}_1-\textbf{v}_2-\textbf{v}_3=&(-a+1,-a+1,-a,a+1)
\end{array}
\]

with $\textbf{v}_1, \textbf{v}_2, \textbf{v}_3, \textbf{v}_4$ as defined for $k \equiv 1$ (mod 2), so that the 16 vectors $\pm \textbf{v}_i$ provide one element of $\textbf{L}_k$ lying strictly within each of the orthants of $\Z^4$. In this case all the coordinates of these vectors have absolute value at most $a+1$.

We must show that each $\textbf{x}\in \Z^4$ is in $S_k+L_k$, which means that for any $\textbf{x} \in \Z^4$ we need to find a $\textbf{w} \in L_k$ such that $\textbf{x}-\textbf{w} \in S_k$. However $\textbf{x}-\textbf{w} \in S_k$ if and only if $\delta (\textbf{x},\textbf{w})\leq k$, where $\delta$ is the $l^1$ metric on $\Z^4$.
If $\textbf{x}, \textbf{y}, \textbf{z} \in \Z^4$ and each coordinate of $\textbf{y}$ lies between the corresponding coordinate of $\textbf{x}$ and $\textbf{z}$ or is equal to one of them, then $\delta (\textbf{x},\textbf{y})+\delta(\textbf{y},\textbf{z})=\delta(\textbf{x},\textbf{z})$. In such a case we say that  \textquotedblleft $\textbf{y}$ lies between $\textbf{x}$ and $\textbf{z}$\textquotedblright.

For any $\textbf{x}\in \Z^4$, we reduce $\textbf{x}$ by adding appropriate elements of $L_k$ until the resulting vector lies within $l^1$-distance $k$ of $\textbf{0}$ or some other element of $L_k$.
The first stage is to reduce $\textbf{x}$ to a vector whose coordinates all have absolute value at most $a+1$. If $\textbf{x}$ has a coordinate with absolute value above $a+1$, then let $\textbf{v}$ be one of the vectors $\pm \textbf{v}_i(1 \leq i \leq 8)$ such that the coordinates of  $\textbf{v}$ have the same sign as the corresponding coordinates of $\textbf{x}$. If a coordinate of $\textbf{x}$ is 0 then either sign is allowed for $\textbf{v}$ as long as the corresponding coordinate of $\textbf{v}$ has absolute value $\leq a+1$. So in the case $k \equiv 0 \pmod 2$ if the $\textbf{e}_3$ coordinate of $\textbf{x}$ is 0 then we  avoid $\textbf{v}_7$ and take $\textbf{v}_5$ instead. Also if the $\textbf{e}_4$ coordinate of $\textbf{x}$ is 0 (or both $\textbf{e}_3$ and $\textbf{e}_4$ coordinates are 0) then instead of $\textbf{v}_5$ we take $\textbf{v}_1$.

Now consider $\textbf{x}'=\textbf{x}-\textbf{v}$. If a coordinate of $\textbf{x}$ has absolute value $s, 1\leq s\leq a+1$, then the corresponding coordinate of $\textbf{x}'$ will have absolute value $s'\leq a+1$ because of the sign matching and the fact that the coordinates of $\textbf{v}$ have absolute value $\leq a+2$. If a coordinate of $\textbf{x}$ has absolute value $s=0$, then as indicated above, the corresponding value of $\textbf{x}'$ will have absolute value $s' \leq a+1$ because $\textbf{v}$ is chosen such that the corresponding coordinate has absolute value $\leq a+1$. If a coordinate of $\textbf{x}$ has absolute value $s>a+1$, then the corresponding coordinate of $\textbf{x}'$ will be strictly smaller in absolute value. Therefore repeating this procedure will result in a vector whose coordinates all have absolute value at most $a+1$.

If the resulting vector $\textbf{x}'$ lies between $\textbf{0}$ and $\textbf{v}$, where $\textbf{v}=\pm\textbf{v}_i$ for some $i$, then we have $\delta(\textbf{0},\textbf{x}')+\delta(\textbf{x}',\textbf{v})=\delta(\textbf{0},\textbf{v})$. For $k\equiv 0 \pmod 2$ all of the vectors $\textbf{v}$ satisfy $\delta(\textbf{0},\textbf{v})=4a+1$, and for $k\equiv 1 \pmod 2$ they all satisfy $\delta(\textbf{0},\textbf{v})=4a-1$. So in either case we have $\delta(\textbf{0},\textbf{v})=2k+1$. Since $\delta(\textbf{0},\textbf{x}')$ and $\delta(\textbf{x}',\textbf{v})$ are both non-negative integers, one of them must be at most $k$, so that $\textbf{x}' \in S_k+L_k$. Hence we also have $\textbf{x} \in S_k+L_k$ as required.

Now we are left with the case where the absolute value of each coordinate of the reduced $\textbf{x}$ is at most $a+1$, and $\textbf{x}$ is in the orthant of $\textbf{v}$, where $\textbf{v} = \pm \textbf{v}_i$ for some $i \leq 8$ but does not lie between $\textbf{0}$ and $\textbf{v}$.
Since $L_k$ is centrosymmetric we only need to consider the eight orthants containing $\textbf{v}_1, ..., \textbf{v}_8$.
For both cases $k \equiv0$ and $k\equiv1 \pmod 2$ the exceptions need to be considered for each orthant in turn. We first consider all eight orthants for the case $ k \equiv 0 \pmod 2$ and then the same for $k \equiv 1 \pmod 2$.

\section*{Orthant of $\textbf{v}_1$, $k \equiv 0$ (mod 2) }
Suppose that $k \equiv 0 \pmod 2$ and $\textbf{x}$ lies within the orthant of $\textbf{v}_1$, but not between $\textbf{0}$ and $\textbf{v}_1$. Then as $\textbf{v}_1 = (-a-1,a+1,a,-a+1)$, the third coordinate of $\textbf{x}$ is equal to $a+1$ or the fourth coordinate equals $-a$ or $-a-1$. We distinguish three cases.

Case 1: $\textbf{x}=(-r,s,a+1,-u)$ where $0\leq r, s\leq a+1$ and $a \leq u \leq a+1$.
Let $\textbf{x}' = \textbf{x} - \textbf{v}_1 = (a+1-r, s-a-1, 1, a-1-u)$, which lies between $\textbf{0}$ and $-\textbf{v}_7$ unless $r \leq 1$ or $s \leq 1$. Let $\textbf{x}'' = \textbf{x}'+\textbf{v}_7=(2-r,s-2,-a-1,2a-u)$.
If $r \leq 1$ and $s \leq 1$ then $\textbf{x}''$ lies between $\textbf{0}$ and $-\textbf{v}_1$ unless $u=a$, in which case let $\textbf{x}''' = \textbf{x}''+ \textbf{v}_1 = (1-a-r, a-1+s, -1, a+1-u)$ which lies between $\textbf{0}$ and $\textbf{v}_7$.
If $r \leq 1$ and $s \geq 2$ then $\textbf{x}''$ lies between $\textbf{0}$ and $-\textbf{v}_3$.
If $r \geq 2$ and $s \leq 1$ then $\textbf{x}''$ lies between $\textbf{0}$ and $-\textbf{v}_2$.

Case 2: $\textbf{x}=(-r,s,a+1,-u)$ where $0\leq r, s\leq a+1$ and $0 \leq u \leq a-1$.
Let $\textbf{x}' = \textbf{x} - \textbf{v}_1 = (a+1-r, s-a-1, 1, a-1-u)$, which lies between $\textbf{0}$ and $-\textbf{v}_6$ unless $r = 0$ or $s = 0$.
Let $\textbf{x}'' = \textbf{x}'+\textbf{v}_6=(1-r,s-1,-a,-u-1)$. If $r = 0$ and $s = 0$ then $\textbf{x}''$ lies between $\textbf{0}$ and $-\textbf{v}_5$.
If $r = 0$ and $s \geq 1$ then $\textbf{x}''$ lies between $\textbf{0}$ and $-\textbf{v}_4$.
If $r \geq 1$ and $s = 0$ then $\textbf{x}''$ lies between $\textbf{0}$ and $-\textbf{v}_8$.

Case 3: $\textbf{x}=(-r,s,t,-u)$ where $0\leq r, s\leq a+1$ and $0 \leq t \leq a$ and $a \leq u \leq a+1$.
Let $\textbf{x}' = \textbf{x} - \textbf{v}_1 = (a+1-r, s-a-1, t-a, a-1-u)$, which lies between $\textbf{0}$ and $-\textbf{v}_5$ unless $r = 0$ or $s = 0$ or $t=0$. If $r=0$ and $s=0$, then $\textbf{x}$ lies between $\textbf{0}$ and $-\textbf{v}_7$.
Let $\textbf{x}'' = \textbf{x}'+\textbf{v}_5=(1-r,s-1,t-1,2a+1-u)$. If $r = 0, s\geq 1$ and $t \geq 1$ then $\textbf{x}''$ lies between $\textbf{0}$ and $\textbf{v}_8$.
Let $\textbf{x}'''=\textbf{x}+\textbf{v}_4 = (-a-r,s-a,a+t,a+1-u)$.
If $r = 0$ and $s \geq 1$ and $t=0$, then $\textbf{x}'''$ lies between $\textbf{0}$ and $\textbf{v}_4$ unless $s=a+1$, in which case if $u=a$ then $\textbf{x}$ lies between $\textbf{0}$ and $\textbf{v}_2$, and if $u=a+1$ then $\textbf{x}'''$ lies between $\textbf{0}$ and $\textbf{v}_4$.
Let $\textbf{x}'''' = \textbf{x}-\textbf{v}_3=(a+1-r,a-1+s,t-a-1,a-u)$.
If $r \geq 1, s= 0$ and $t \geq 1$ then $\textbf{x}''''$ lies between $\textbf{0}$ and $-\textbf{v}_4$. If $r \geq 1, s= 0$ and $t = 0$ then $\textbf{x}''''$ lies between $\textbf{0}$ and $-\textbf{v}_3$ if $u=a$, and between $\textbf{0}$ and $\textbf{v}_6$ if $u=a+1$.
If $r \geq 1, s \geq 1$ and $t = 0$ then $\textbf{x}''$ lies between $\textbf{0}$ and $\textbf{v}_7$ unless $r=a+1$ or $s=a+1$.
If $r = a+1$, $s \geq 1$ and $t =0$ then $\textbf{x}'$ lies between $\textbf{0}$ and $-\textbf{v}_8$.
If $r \geq 1$, $s=a+1$ and $t =0$ then $\textbf{x}'$ lies between $\textbf{0}$ and $-\textbf{v}_4$.

This completes the cases for the orthant of $\textbf{v}_1$ for $k \equiv 0 \pmod 2$.

\section*{Orthant of $\textbf{v}_2$, $k \equiv 0$ (mod 2) }
Now suppose that $\textbf{x}$ lies in the orthant of $\textbf{v}_2$ but not between $\textbf{0}$ and $\textbf{v}_2$. Then the first coordinate of $\textbf{x}$ is equal to $a$ or $a+1$, or the fourth coordinate equals $-a-1$. We distinguish three cases.

Case 1: $\textbf{x}=(r,s,t,-a-1)$ where $a \leq r \leq a+1$ and $0 \leq s,t \leq a+1$.
Let $\textbf{x}' = \textbf{x} - \textbf{v}_2 = (r-a+1, s-a-1, t-a-1, -1)$, which lies between $\textbf{0}$ and $-\textbf{v}_5$ unless $s=0$ or $t \leq 1$, in which case let $\textbf{x}'' = \textbf{x}'+\textbf{v}_5=(r-2a+1, s-1, t-2, a+1)$.
If $s=0$ and $t \leq 1$ then let $\textbf{x}'''= \textbf{x}+ \textbf{v}_5 = (r-a,a,t+a-1,1)$ which lies between $\textbf{0}$ and $\textbf{v}_8$.
If $s=0$ and $t \geq 2$ then let $\textbf{x}''''= \textbf{x}'''- \textbf{v}_8 = (r-2a,0,t-1,-a)$ which lies between $\textbf{0}$ and $\textbf{v}_3$.
If $s \geq 1$ and $t \leq 1$ then $\textbf{x}''$ lies between $\textbf{0}$ and $\textbf{v}_7$ unless $s=a+1$, in which case let $\textbf{x}^v= \textbf{x}''- \textbf{v}_7 = (r-a,1,a+t,0)$ which lies between $\textbf{0}$ and $\textbf{v}_2$ .

Case 2: $\textbf{x}=(r,s,t,-u)$ where $a\leq r \leq a+1$, $0 \leq s,t \leq a+1$ and $0 \leq u \leq a$.
Let $\textbf{x}' = \textbf{x} - \textbf{v}_2 = (r-a+1,s-a-1, t-a-1, a-u)$, which lies between $\textbf{0}$ and $-\textbf{v}_1$ unless $t = 0$ or $u = 0$.
If $t = 0$ and $u = 0$ then $\textbf{x}$ lies between $\textbf{0}$ and $ -\textbf{v}_3$ unless $a \leq s \leq a+1$.
If $r=a+1$, $a \leq s \leq a+1$, $t=0$ and $u = 0$ then let $\textbf{x}''$ = $\textbf{x} + \textbf{v}_3 = (r-a-1, s-a+1,t+a+1,-u-a)$ which lies between $\textbf{0}$ and $\textbf{v}_2$.
If $a \leq r \leq a+1$, $s=a+1$, $t=0$ and $u=0$ then let $\textbf{x}'' = \textbf{x}-\textbf{v}_2=(r-a+1,s-a-1,t-a-1,-u+a)$ which lies between $\textbf{0}$ and $- \textbf{v}_3$.
If $r=a$, $s=a$, $t=0$ and $u=0$ then $\textbf{x}$ lies between $\textbf{0}$ and $\textbf{v}_8$.
Now let $\textbf{x}'''= \textbf{x}'+\textbf{v}_1=(r-2a,s,t-1,1-u)$.
If $t=0$ and $1 \leq u \leq a$ then $\textbf{x}'''$ lies between $\textbf{0}$ and $\textbf{v}_6$ unless $s = a+1$, in which case let $\textbf{x}'''' = \textbf{x}''' - \textbf{v}_6 = (r-a,s-a,t+a,a+1-u)$ which lies between $\textbf{0}$ and $\textbf{v}_8$.
If $1 \leq t \leq a+1$ and $u=0$ then $\textbf{x}'''$ lies between $\textbf{0}$ and $\textbf{v}_5$ unless $s = a+1$ or $t=a+1$, in which case let $\textbf{x}'''' = \textbf{x}'''-\textbf{v}_5=(r-a,s-a,t-a,-a-1-u)$. 
If $s=a+1$ and $t=a+1$ then $ \textbf{x}'$ lies between $\textbf{0}$ and $\textbf{v}_8$.
If $s=a+1$ and $1 \leq t \leq a$ then $ \textbf{x}''''$ lies between $\textbf{0}$ and $- \textbf{v}_4$. If $0 \leq s \leq a$ and $t=a+1$ then $ \textbf{x}''''$ lies between $\textbf{0}$ and $- \textbf{v}_7$ unless $s=0$, in which case $\textbf{x}'''$ lies between $\textbf{0}$ and $\textbf{v}_4$.

Case 3: $\textbf{x}=(r,s,t,-a-1)$ where $0 \leq r \leq a-1$ and $0 \leq s,t \leq a+1$.
Let $\textbf{x}' = \textbf{x} - \textbf{v}_2 = (r-a+1, s-a-1, t-a-1, -1)$, which lies between $\textbf{0}$ and $-\textbf{v}_8$ unless $s = 0$ or $t = 0$.
If $s=0$ then $\textbf{x}$ lies between $\textbf{0}$ and $-\textbf{v}_7$.
If $t=0$ then $\textbf{x}$ lies between $\textbf{0}$ and $-\textbf{v}_4$ unless $s=a+1$, in which case let $\textbf{x}'' = \textbf{x}+\textbf{v}_4=(r-a,1,a,0)$ which lies between $\textbf{0}$ and $\textbf{v}_1$. \newline
This completes the cases for the orthant of $\textbf{v}_2$.

\section*{Orthant of $\textbf{v}_3$, $k \equiv 0$ (mod 2) }
Now suppose that $\textbf{x}$ lies in the orthant of $\textbf{v}_3$ but not between $\textbf{0}$ and $\textbf{v}_3$. Then the second coordinate of $\textbf{x}$ is equal to $-a$ or $-a-1$, or the fourth coordinate equals $-a-1$. We distinguish three cases.

Case 1: $\textbf{x}=(-r,-s,t,-a-1)$ where $0 \leq r,t \leq a+1$ and $a \leq s \leq a+1$.
Let $\textbf{x}' = \textbf{x} - \textbf{v}_3 = (a+1-r, a-1-s, t-a-1, -1)$, which lies between $\textbf{0}$ and $-\textbf{v}_5$ unless $r = 0$ or $t \leq 1$, in which case let $\textbf{x}'' = \textbf{x}'+\textbf{v}_5=(1-r, 2a-1-s, t-2, a+1)$.
If $r=0$ and $t \geq 2$ then $\textbf{x}''$ which lies between $\textbf{0}$ and $\textbf{v}_8$.
If $r=0$ and $t \leq 1$ then let $\textbf{x}'''= \textbf{x} + \textbf{v}_5 = (-a,a-s,a-1+t,1)$ which lies between $\textbf{0}$ and $\textbf{v}_4$.
If $r \geq 1$ and $t \leq 1$ then $\textbf{x}''$ lies between $\textbf{0}$ and $\textbf{v}_7$ unless $r=a+1$, in which case let $\textbf{x}''''= \textbf{x}''- \textbf{v}_7 = (-1,a-s,a+t,0)$ which lies between $\textbf{0}$ and $\textbf{v}_3$.

Case 2: $\textbf{x}=(-r,-s,t,-a-1)$ where $0 \leq r,t \leq a+1$ and $0 \leq s \leq a-1$.
Let $\textbf{x}' = \textbf{x} - \textbf{v}_3 = (a+1-r, a-1-s, t-a-1, -1)$, which lies between $\textbf{0}$ and $-\textbf{v}_4$ unless $r = 0$ or $t = 0$. Let $\textbf{x}'' = \textbf{x}'+\textbf{v}_4=(1-r,-1-s,t-1,a)$.
If $r = 0$ and $t \geq 1$ then $\textbf{x}''$ lies between $\textbf{0}$ and $-\textbf{v}_6$.
If $r \geq 1$ and $t = 0$ then $\textbf{x}''$ lies between $\textbf{0}$ and $-\textbf{v}_2$ unless $r=a+1$, in which case let $\textbf{x}''' = \textbf{x}''+\textbf{v}_2=(-1,a-s,a,0)$ which lies between $\textbf{0}$ and $\textbf{v}_1$.
If $r = 0$ and $t =0$, then let $\textbf{x}''' = \textbf{x}''+\textbf{v}_1=(-a,a-s,a-1,1)$ which lies between $\textbf{0}$ and $\textbf{v}_5$.

Case 3: $\textbf{x}=(-r,-s,t,-u)$ where $0 \leq r,t \leq a+1$, $a \leq s \leq a+1$ and $0 \leq u \leq a$.
Let $\textbf{x}' = \textbf{x} - \textbf{v}_3 = (a+1-r, a-1-s, t-a-1, a-u)$, which lies between $\textbf{0}$ and $-\textbf{v}_1$ unless $t = 0$ or $u = 0$.
Let $\textbf{x}'' = \textbf{x} + \textbf{v}_8 = (a-r, a-s, a+t, a+1-u)$. 
If $t=0$ and $u=0$ then $\textbf{x}''$ lies between $\textbf{0}$ and $\textbf{v}_4$ unless $r \leq a-1$, in which case let $\textbf{x}''' = \textbf{x}+\textbf{v}_2=(a-1-r,a+1-s,a+1+t,-a-u)$ which lies between $\textbf{0}$ and $\textbf{v}_2$.
If $t=0$ and $u \geq 1$ then $\textbf{x}''$ lies between $\textbf{0}$ and $-\textbf{v}_6$ unless $r = a+1$, in which case $\textbf{x}''$ lies between $\textbf{0}$ and $\textbf{v}_4$.
Let $\textbf{x}'''' = \textbf{x}-\textbf{v}_4=(a-r,a-s,t-a,-a-1-u)$. If $t \geq 1, u=0$ and $r \leq a$ then $\textbf{x}''''$ lies between $\textbf{0}$ and $-\textbf{v}_5$ unless $t = a+1$. If $t = a+1, u=0$ and $r \leq a$ then $\textbf{x}''''$ lies between $\textbf{0}$ and $-\textbf{v}_7$ unless $r = 0$ in which case $\textbf{x}''$ lies between $\textbf{0}$ and $\textbf{v}_8$. If $t \geq 1, u=0$ and $r = a+1$ then $\textbf{x}'$ lies between $\textbf{0}$ and  $-\textbf{v}_2$. \newline
This completes the cases for the orthant of $\textbf{v}_3$.

\section*{Orthant of $\textbf{v}_4$, $k \equiv 0$ (mod 2) }
Now suppose $\textbf{x}$ lies in the orthant of $\textbf{v}_4$ but not between $\textbf{0}$ and $\textbf{v}_4$. Then the first coordinate of $\textbf{x}$ is equal to $-a-1$ or the second coordinate is equal to $-a-1$, or the third equals $a+1$. We distinguish seven cases.

Case 1: $\textbf{x}=(-a-1,-a-1,a+1,u)$ where $0 \leq u \leq a+1$.
Let $\textbf{x}' = \textbf{x} - \textbf{v}_4 = (-1, -1, 1,u-a-1)$, which lies between $\textbf{0}$ and $\textbf{v}_4$ if $u=a+1$ and between $\textbf{0}$ and $\textbf{v}_3$ if $u\leq a$.

Case 2: $\textbf{x}=(-a-1,-a-1,t,u)$ where $0 \leq t \leq a$ and $0 \leq u \leq a+1$.
Let $\textbf{x}' = \textbf{x} - \textbf{v}_4 = (-1, -1, t-a,u-a-1)$, which lies between $\textbf{0}$ and $-\textbf{v}_8$.

Case 3: $\textbf{x}=(-a-1,-s,a+1,u)$ where $0 \leq s \leq a$ and $0 \leq u \leq a+1$.
Let $\textbf{x}' = \textbf{x} - \textbf{v}_4 = (-1, -a-s, 1, u-a-1)$, which lies between $\textbf{0}$ and $\textbf{v}_1$ unless $u \geq a$, in which case let $\textbf{x}''= \textbf{x}'- \textbf{v}_1 = (a,-s-1,-a+1,u-2)$ which lies between $\textbf{0}$ and $-\textbf{v}_1$.

Case 4: $\textbf{x}=(-r,-a-1,a+1,u)$ where $0 \leq r \leq a$ and $0 \leq u \leq a+1$.
Let $\textbf{x}' = \textbf{x} - \textbf{v}_4 = (a-r, -1, 1, u-a-1)$, which lies between $\textbf{0}$ and $-\textbf{v}_7$ unless $r = 0$, in which case let $\textbf{x}''= \textbf{x}'+ \textbf{v}_7 = (1,a-2,-a-1,u)$ which lies between $\textbf{0}$ and $-\textbf{v}_3$ unless $u = a+1$, in which case let $\textbf{x}'''= \textbf{x}''+ \textbf{v}_3 = (-a,-1,0,1)$ which lies between $\textbf{0}$ and $\textbf{v}_4$.

Case 5: $\textbf{x}=(-r,-s,a+1,u)$ where $0 \leq r,s \leq a$ and $0 \leq u \leq a+1$.
Let $\textbf{x}' = \textbf{x} - \textbf{v}_4 = (a-r, a-s, 1,u-a-1)$, which lies between $\textbf{0}$ and $\textbf{v}_2$ unless $r = 0$ or $u = 0$ in which case let $\textbf{x}'' = \textbf{x}'-\textbf{v}_2=(1-r,-s-1,-a,u-1)$.
If $r \geq 1$ and $u = 0$ then $\textbf{x}''$ lies between $\textbf{0}$ and $-\textbf{v}_8$ unless $s=a$, in which case let $\textbf{x}''' = \textbf{x}''+\textbf{v}_8=(a+1-r,-1,0,a)$ which lies between $\textbf{0}$ and $-\textbf{v}_6$.
If $r = 0$ then $\textbf{x}''$ lies between $\textbf{0}$ and $-\textbf{v}_1$ unless $u=0$ or $u=a+1$.
If $r = 0$ and $u =0$, then let $\textbf{x}'''' = \textbf{x}''+\textbf{v}_5=(1-a,a-1-s,-1,a+1)$ which lies between $\textbf{0}$ and $\textbf{v}_7$ unless $s=a$, in which case let $\textbf{x}^v = \textbf{x}''''+\textbf{v}_2=(0,a,a,1)$ which lies between $\textbf{0}$ and $\textbf{v}_8$.
If $r = 0$ and $u =a+1$, then let $\textbf{x}^{v\imath} = \textbf{x}''+\textbf{v}_1=(-a,a-s,0,1)$ which lies between $\textbf{0}$ and $\textbf{v}_5$.

Case 6: $\textbf{x}=(-r,-a-1,t,u)$ where $0 \leq r,t \leq a$ and $0 \leq u \leq a+1$.
Let $\textbf{x}' = \textbf{x} - \textbf{v}_4 = (a-r, -1, t-a,u-a-1)$, which lies between $\textbf{0}$ and $-\textbf{v}_5$ unless $t = 0$, in which case let $\textbf{x}'' = \textbf{x}'+\textbf{v}_5=(-r, a-1,-1,u+1)$ which lies between $\textbf{0}$ and $\textbf{v}_7$ unless $r = a$ or $u=a+1$.
If $t = 0$ and $r = a$ then let $\textbf{x}''' = \textbf{x}''-\textbf{v}_7=(-1,0,a+1,u-a)$ which lies between $\textbf{0}$ and $\textbf{v}_3$ unless $u=a+1$.
If $t = 0$ and $u=a+1$ then let $\textbf{x}''' = \textbf{x}''-\textbf{v}_7=(a-1-r,0,a+1,1)$ which lies between $\textbf{0}$ and $-\textbf{v}_6$ unless $r=a$.
If $t = 0, r=a$ and $u =a+1$, then $\textbf{x}''' = (-1,0,a+1,1)$. 
Let $\textbf{x}'''' = \textbf{x}'''-\textbf{v}_4=(a-1,a,1,-a)$ which lies between $\textbf{0}$ and $\textbf{v}_2$.

Case 7: $\textbf{x}=(-a-1,-s,t,u)$ where $0 \leq s,t \leq a$ and $0 \leq u \leq a+1$.
Let $\textbf{x}' = \textbf{x} - \textbf{v}_4 = (-1, a-s, t-a,u-a-1)$, which lies between $\textbf{0}$ and $\textbf{v}_6$ unless $u = 0$, in which case let $\textbf{x}'' = \textbf{x}'-\textbf{v}_6=(a-1, -s, t+1, -1)$ which lies between $\textbf{0}$ and $-\textbf{v}_7$ unless $s = a$, in which case let $\textbf{x}''' = \textbf{x}''+\textbf{v}_7=(0,-1,t-a-1,a)$ which lies between $\textbf{0}$ and $-\textbf{v}_2$. \newline 
This completes the cases for the orthant of $\textbf{v}_4$.

\section*{Orthant of $\textbf{v}_5$, $k \equiv 0$ (mod 2) }
Now suppose $\textbf{x}$ lies in the orthant of $\textbf{v}_5$ but not between $\textbf{0}$ and $\textbf{v}_5$. Then the first coordinate of $\textbf{x}$ is equal to $-a-1$ or the second coordinate is equal to $a+1$, or the third equals $a$ or $a+1$. We distinguish seven cases.

Case 1: $\textbf{x}=(-a-1,a+1,t,u)$ where $a \leq t \leq a+1$ and $0 \leq u \leq a+1$.
Let $\textbf{x}' = \textbf{x} - \textbf{v}_5 = (-1, 1, t-a+1,u-a-2)$, which lies between $\textbf{0}$ and $\textbf{v}_1$ unless $u \leq 2$, in which case let $\textbf{x}''= \textbf{x}'- \textbf{v}_1 = (a,-a,t-2a+1,u-3)$ which lies between $\textbf{0}$ and $-\textbf{v}_5$.

Case 2: $\textbf{x}=(-a-1,a+1,t,u)$ where $0 \leq t \leq a-1$ and $0 \leq u \leq a+1$.
Let $\textbf{x}' = \textbf{x} - \textbf{v}_5 = (-1, 1, t-a+1,u-a-2)$, which lies between $\textbf{0}$ and $\textbf{v}_6$ unless $u \leq 1$, in which case let $\textbf{x}''= \textbf{x}'- \textbf{v}_6 = (a-1,-a+1,t+2,u-2)$ which lies between $\textbf{0}$ and $-\textbf{v}_7$.

Case 3: $\textbf{x}=(-a-1,s,t,u)$ where $0 \leq s \leq a$, $a \leq t \leq a+1$ and $0 \leq u \leq a+1$.
Let $\textbf{x}' = \textbf{x} - \textbf{v}_5 = (-1, s-a, t-a+1,u-a-2)$, which lies between $\textbf{0}$ and $\textbf{v}_3$ unless $s = 0$ or $u \leq 1$, in which case let $\textbf{x}''= \textbf{x}'- \textbf{v}_3 = (a,s-1,t-2a,u-2)$.
If $s = 0$ and $u \leq 1$ then $\textbf{x}''$ lies between $\textbf{0}$ and $-\textbf{v}_5$ unless $t=a$, in which case let $\textbf{x}''' = \textbf{x}''+\textbf{v}_5=(0,a-1,-1,a+u)$ which lies between $\textbf{0}$ and $\textbf{v}_7$.
If $s = 0$ and $u \geq 2$ then let $\textbf{x}''' = \textbf{x}''+\textbf{v}_1=(-1,a,t-a,u-a-1)$ which lies between $\textbf{0}$ and $\textbf{v}_1$.
If $s \geq 1$ and $u \leq 1$ then $\textbf{x}''$ lies between $\textbf{0}$ and $-\textbf{v}_4$.

Case 4: $\textbf{x}=(-r,a+1,t,u)$ where $0 \leq r \leq a$, $a \leq t \leq a+1$ and $0 \leq u \leq a+1$.
Let $\textbf{x}' = \textbf{x} - \textbf{v}_5 = (a-r, 1, t-a+1, u-a-2)$, which lies between $\textbf{0}$ and $\textbf{v}_2$ unless $r = 0$ or $u \geq 2$, in which case let $\textbf{x}''= \textbf{x}'- \textbf{v}_2 = (1-r,-a,t-2a,u-2)$. 
If $r = 0$ and $u \geq 2$ then $\textbf{x}''$ lies between $\textbf{0}$ and $-\textbf{v}_1$.
If $r = 0$ and $u \leq 1$ then $\textbf{x}''$ lies between $\textbf{0}$ and $-\textbf{v}_5$ unless $t=a$, in which case let $\textbf{x}''' = \textbf{x}''+\textbf{v}_5=(1-a,0,-1,a+u)$ which lies between $\textbf{0}$ and $\textbf{v}_7$.
If $r \geq 1$ and $u \geq 2$ then $\textbf{x}''$ lies between $\textbf{0}$ and $-\textbf{v}_2$.

Case 5: $\textbf{x}=(-r,s,t,u)$ where $0 \leq r,s \leq a$, $a \leq t \leq a+1$ and $0 \leq u \leq a+1$.
Let $\textbf{x}' = \textbf{x} - \textbf{v}_5 = (a-r, s-a, t-a+1,u-a-2)$, which lies between $\textbf{0}$ and $-\textbf{v}_7$ unless $r=0$ or $s=0$ or $u=0$.
If $r = 0$ then $\textbf{x}$ lies between $\textbf{0}$ and $\textbf{v}_8$ unless $t=a+1$, in which case let $\textbf{x}' = \textbf{x}-\textbf{v}_8=(-a,s-a,1,u-a-1)$ which lies between $\textbf{0}$ and $\textbf{v}_3$ unless $u = 0$.
If $r = 0, t=a+1$ and $u=0$ then $\textbf{x}$ lies between $\textbf{0}$ and $\textbf{v}_2$.
If $r \geq 1$ and $s =0$ then $\textbf{x}$ lies between $\textbf{0}$ and $\textbf{v}_4$ unless $t=a+1$, in which case let $\textbf{x}' = \textbf{x}-\textbf{v}_4=(a-r,a,1,u-a-1)$ which lies between $\textbf{0}$ and $\textbf{v}_2$ unless $u=0$.
If $r \geq 1, s=0$ and $u = 0$, then let $\textbf{x}'' = \textbf{x}'-\textbf{v}_2=(1-r,-1,-a,-1)$ which lies between $\textbf{0}$ and $-\textbf{v}_8$.
If $r \geq 1, s \geq 1$ and $u =0$, then $\textbf{x}$ lies between $\textbf{0}$ and $\textbf{v}_1$ unless $t=a+1$, in which case let $\textbf{x}' = \textbf{x}-\textbf{v}_1=(a+1-r,s-a-1,1,a-1)$ which lies between $\textbf{0}$ and $-\textbf{v}_7$.

Case 6: $\textbf{x}=(-r,a+1,t,u)$ where $0 \leq r \leq a$, $0 \leq t \leq a-1$ and $0 \leq u \leq a+1$.
Let $\textbf{x}' = \textbf{x} - \textbf{v}_5 = (a-r, 1, t-a+1,u-a-2)$, which lies between $\textbf{0}$ and $-\textbf{v}_4$ unless $u = 0$ in which case $\textbf{x}$ lies between $\textbf{0}$ and $\textbf{v}_1$.

Case 7: $\textbf{x}=(-a-1,s,t,u)$ where $0 \leq s \leq a$, $0 \leq t \leq a-1$ and $0 \leq u \leq a+1$.
Let $\textbf{x}' = \textbf{x} - \textbf{v}_5 = (-1, s-a, t-a+1,u-a-2)$, which lies between $\textbf{0}$ and $-\textbf{v}_8$ unless $u = 0$ in which case $\textbf{x}$ lies between $\textbf{0}$ and $\textbf{v}_1$. \newline
This completes the cases for the orthant of $\textbf{v}_5$.

\section*{Orthant of $\textbf{v}_6$, $k \equiv 0$ (mod 2) }
Now suppose $\textbf{x}$ lies in the orthant of $\textbf{v}_6$ but not between $\textbf{0}$ and $\textbf{v}_6$. Then the first coordinate of $\textbf{x}$ is equal to $-a-1$ or the second coordinate is equal to $a+1$, or the fourth equals $-a-1$. We distinguish seven cases.

Case 1: $\textbf{x}=(-a-1,a+1,-t,-a-1)$ where $0 \leq t \leq a+1$.
Let $\textbf{x}' = \textbf{x} - \textbf{v}_6 = (-1, 1, a+1-t,-1)$, which lies between $\textbf{0}$ and $\textbf{v}_1$ unless $t = 0$, in which case let $\textbf{x}''= \textbf{x}'- \textbf{v}_1 = (a,-a,1,a-2)$ which lies between $\textbf{0}$ and $-\textbf{v}_6$.

Case 2: $\textbf{x}=(-a-1,a+1,-t,-u)$ where $0 \leq t \leq a+1$ and $0 \leq u \leq a$.
Let $\textbf{x}' = \textbf{x} - \textbf{v}_6 = (-1, 1, a+1-t,a-u)$, which lies between $\textbf{0}$ and $\textbf{v}_5$ unless $t \leq 1$, in which case let $\textbf{x}''= \textbf{x}'- \textbf{v}_5 = (a-1,1-a,2-t,-u-2)$ which lies between $\textbf{0}$ and $-\textbf{v}_7$ unless $u = a$.
If $t = 1$ and $u = a$ then $\textbf{x}'$ lies between $\textbf{0}$ and $\textbf{v}_1$
If $t = 0$ and $u = a$ then let $\textbf{x}'' = \textbf{x}'-\textbf{v}_1=(a,-a,1,a+1)$ and $\textbf{x}''' = \textbf{x}''+\textbf{v}_6=(0,0,-a,1)$
which lies between $\textbf{0}$ and $\textbf{v}_7$.

Case 3: $\textbf{x}=(-a-1,s,-t,-a-1)$ where $0 \leq s \leq a$ and $0 \leq t \leq a+1$.
Let $\textbf{x}' = \textbf{x} - \textbf{v}_6 = (-1, s-a, a+1-t, -1)$, which lies between $\textbf{0}$ and $\textbf{v}_3$ unless $s = a$ in which case let $\textbf{x}''= \textbf{x}'- \textbf{v}_3 = (a,a-1,-t,a-1)$ which lies between $\textbf{0}$ and $-\textbf{v}_3$.

Case 4: $\textbf{x}=(-r,a+1,-t,-a-1)$ where $0 \leq r \leq a$ and $0 \leq t \leq a+1$.
Let $\textbf{x}' = \textbf{x} - \textbf{v}_6 = (a-r, 1, a+1-t,-1)$, which lies between $\textbf{0}$ and $\textbf{v}_2$ unless $r = 0$, in which case let $\textbf{x}''= \textbf{x}'- \textbf{v}_2 = (1,-a,-t,-a-1)$ which lies between $\textbf{0}$ and $-\textbf{v}_1$ unless $t = a+1$ in which case $\textbf{x}' = (a,1,0,-1)$ which lies between $\textbf{0}$ and $-\textbf{v}_4$

Case 5: $\textbf{x}=(-r,s,-t,-a-1)$ where $0 \leq r,s \leq a$ and $0 \leq t \leq a+1$.
Let $\textbf{x}' = \textbf{x} - \textbf{v}_6 = (a-r, s-a, a+1-t, -1)$, which lies between $\textbf{0}$ and $-\textbf{v}_7$ unless $r = 0$ or $s = 0$, in which case let $\textbf{x}''= \textbf{x}'+ \textbf{v}_7 = (1-r,s-1,-t-1,a)$. 
If $r = 0$ and $s \geq 1$ then $\textbf{x}''$ lies between $\textbf{0}$ and $-\textbf{v}_3$ unless $t=a+1$, in which case let $\textbf{x}''' = \textbf{x}''+\textbf{v}_3=(-a,s-a,-1,0)$ which lies between $\textbf{0}$ and $\textbf{v}_6$.
If $r \geq 1$ and $s = 0$ then $\textbf{x}''$ lies between $\textbf{0}$ and $-\textbf{v}_2$ unless $t=a+1$, in which case let $\textbf{x}''' = \textbf{x}''+\textbf{v}_2=(a-r, a,-1,0)$ which lies between $\textbf{0}$ and $-\textbf{v}_4$.
If $r =0$ and $s = 0$ then $\textbf{x}$ lies between $\textbf{0}$ and $-\textbf{v}_8$ unless $t=a+1$, in which case let $\textbf{x}'''' = \textbf{x}+\textbf{v}_8=(a,a,-1,0)$ which lies between $\textbf{0}$ and $-\textbf{v}_4$.

Case 6: $\textbf{x}=(-r,a+1,-t,-u)$ where $0 \leq r,u \leq a$ and $0 \leq t \leq a+1$.
Let $\textbf{x}' = \textbf{x} - \textbf{v}_6 = (a-r, 1, a+1-t,a-u)$, which lies between $\textbf{0}$ and $\textbf{v}_8$ unless $t=0$ in which case let $\textbf{x} = (-r,a+1,0,-u)$ which lies between $\textbf{0}$ and $\textbf{v}_1$ unless $u=a$.
If $t = 0$ and $u=a$ then let $\textbf{x}'' = \textbf{x}-\textbf{v}_1=(a+1-r,0,-a,-1)$ which lies between $\textbf{0}$ and $-\textbf{v}_4$ unless $r = 0$ in which case $\textbf{x} = (0,a+1,0,-a)$ which lies between $\textbf{0}$ and $\textbf{v}_1$.

Case 7: $\textbf{x}=(-a-1,s,-t,-u)$ where $0 \leq s,u \leq a$ and $0 \leq t \leq a+1$.
Let $\textbf{x}' = \textbf{x} - \textbf{v}_6 = (-1, s-a, a+1-t,a-u)$, which lies between $\textbf{0}$ and $\textbf{v}_4$ unless $t = 0$ in which case  $\textbf{x} = (-a-1,s,0,-u)$ which lies between $\textbf{0}$ and $\textbf{v}_1$ unless $u=a$.
If $t = 0$ and $u=a$ then let $\textbf{x}'' = \textbf{x}-\textbf{v}_1=(0,s-a-1,-a,-1)$ which lies between $\textbf{0}$ and $-\textbf{v}_8$ unless $s = 0$ in which case let $\textbf{x}''' = \textbf{x}'' + \textbf{v}_8 = (a,-1,0,a)$ which lies between $\textbf{0}$ and $-\textbf{v}_6$. \newline
This completes the cases for the orthant of $\textbf{v}_6$.

\section*{Orthant of $\textbf{v}_7$, $k \equiv 0$ (mod 2) }
Now suppose $\textbf{x}$ lies in the orthant of $\textbf{v}_7$ but not between $\textbf{0}$ and $\textbf{v}_7$. Then the first coordinate of $\textbf{x}$ is equal to $-a$ or $-a-1$ or the second equals $a$ or $a+1$. We distinguish seven cases.

Case 1: $\textbf{x}=(-r,s,-t,u)$ where $a \leq r, s \leq a+1$, $2 \leq t \leq a+1$ and $0 \leq u \leq a+1$.
Let $\textbf{x}' = \textbf{x} - \textbf{v}_7 = (a-1-r, s-a+1, a+2-t,u-a-1)$, which lies between $\textbf{0}$ and $\textbf{v}_1$ unless $u \leq 1$, in which case let $\textbf{x}''= \textbf{x}- \textbf{v}_6 = (a-r,s-a,a+1-t,a+u)$ which lies between $\textbf{0}$ and $\textbf{v}_5$.

Case 2: $\textbf{x}=(-a,a,-t,u)$ where $0 \leq t \leq 1$ and $0 \leq u \leq a+1$.
If $u = 0$ then $\textbf{x}$ lies between $\textbf{0}$ and $\textbf{v}_6$.
If $u \geq 1$ then let $\textbf{x}' = \textbf{x} - \textbf{v}_7 - \textbf{v}_3 = (a, a, 1-t,u-1)$, which lies between $\textbf{0}$ and $\textbf{v}_8$.

Case 3: $\textbf{x}=(-a-1,a,-t,u)$ where $0 \leq t \leq 1$ and $0 \leq u \leq a+1$.
Let $\textbf{x}' = \textbf{x} - \textbf{v}_1 = (0, -1, -a-t,a-1+u)$.
If $u \leq 1$ then $\textbf{x}'$ lies between $\textbf{0}$ and $-\textbf{v}_2$.
If $u \geq 2$ then let $\textbf{x}'' = \textbf{x}' + \textbf{v}_2 = (a-1, a, 1-t,u-1)$, which lies between $\textbf{0}$ and $\textbf{v}_8$.

Case 4: $\textbf{x}=(-a,a+1,-t,u)$ where $0 \leq t \leq 1$ and $0 \leq u \leq a+1$.
Let $\textbf{x}' = \textbf{x} - \textbf{v}_1 = (1, 0, -a-t,a-1+u)$.
If $u \leq 1$ then $\textbf{x}'$ lies between $\textbf{0}$ and $-\textbf{v}_3$.
If $u \geq 2$ then let $\textbf{x}'' = \textbf{x}' + \textbf{v}_3 = (-a, -a+1, 1-t,u-1)$, which lies between $\textbf{0}$ and $\textbf{v}_4$.

Case 5: $\textbf{x}=(-a-1,a+1,-t,u)$ where $0 \leq t \leq 1$ and $0 \leq u \leq a+1$.
Let $\textbf{x}' = \textbf{x} - \textbf{v}_1 = (0, 0, -a-t,a-1+u)$.
If $u \leq 1$ then $\textbf{x}'$ lies between $\textbf{0}$ and $\textbf{v}_7$.
If $u \geq 2$ then let $\textbf{x}'' = \textbf{x}' - \textbf{v}_7 = (a-1, -a+1, 2-t,u-2)$, which lies between $\textbf{0}$ and $-\textbf{v}_6$.

Case 6: $\textbf{x}=(-r,s,-t,u)$ where $0 \leq r \leq a-1$, $a \leq s \leq a+1$ and $0 \leq t,u \leq a+1$.
Let $\textbf{x}' = \textbf{x} - \textbf{v}_7 = (a-1-r, s-a+1, a+2-t, u-a-1)$, which lies between $\textbf{0}$ and $\textbf{v}_2$ unless $t = 0$ or $u = 0$, in which case let $\textbf{x}''= \textbf{x}'- \textbf{v}_2 = (-r,s-2a,1-t,u-1)$. 
If $t = 0$ and $u \geq 1$ then $\textbf{x}''$ lies between $\textbf{0}$ and $\textbf{v}_4$.
If $t = 0$ and $u = 0$ then $\textbf{x}''$ lies between $\textbf{0}$ and $\textbf{v}_3$ unless $s=a$, in which case let $\textbf{x}''' = \textbf{x}''-\textbf{v}_3=(a+1-r, -1,-a,a-1)$ which lies between $\textbf{0}$ and $-\textbf{v}_1$.
If $t \geq 1$ and $u = 0$ then $\textbf{x}'' = (-r,s-2a,1-t,-1)$ which lies between $\textbf{0}$ and $-\textbf{v}_8$.

Case 7: $\textbf{x}=(-r,s,-t,u)$ where $a \leq r \leq a+1$, $0 \leq s \leq a-1$ and $0 \leq t,u \leq a+1$.
Let $\textbf{x}' = \textbf{x} - \textbf{v}_7 = (a-1-r, s-a+1, a+2-t,u-a-1)$, which lies between $\textbf{0}$ and $\textbf{v}_3$ unless $t=0$ or $u=0$, in which case let $\textbf{x}'' = \textbf{x}' - \textbf{v}_3 = (2a-r,s,1-t,u-1)$.
If $t = 0$ and $u=0$ then $\textbf{x}$ lies between $\textbf{0}$ and $\textbf{v}_1$. 
If $t = 0$ and $u \geq 1$ then $\textbf{x}''$ lies between $\textbf{0}$ and $\textbf{v}_8$. 
If $t \geq 1$ and $u=0$ then $\textbf{x}''$ lies between $\textbf{0}$ and $-\textbf{v}_4$. \newline
This completes the cases for the orthant of $\textbf{v}_7$.

\section*{Orthant of $\textbf{v}_8$, $k \equiv 0$ (mod 2) }
Finally suppose $\textbf{x}$ lies in the orthant of $\textbf{v}_8$ but not between $\textbf{0}$ and $\textbf{v}_8$. Then at least one of the first three coordinate of $\textbf{x}$ is equal to $a+1$. We distinguish seven cases.

Case 1: $\textbf{x}=(a+1,a+1,a+1,u)$ where $0 \leq u \leq a+1$.
Let $\textbf{x}' = \textbf{x} - \textbf{v}_8 = (1, 1, 1,u-a-1)$, which lies between $\textbf{0}$ and $\textbf{v}_2$ unless $u = 0$, in which case let $\textbf{x}''= \textbf{x}' - \textbf{v}_2 = (-a+2,-a,-a,-1)$ which lies between $\textbf{0}$ and $-\textbf{v}_8$.

Case 2: $\textbf{x}=(a+1,a+1,t,u)$ where $0 \leq t \leq a$ and $0 \leq u \leq a+1$.
Let $\textbf{x}' = \textbf{x} - \textbf{v}_8 = (1, 1, t-a,u-a-1)$, which lies between $\textbf{0}$ and $\textbf{v}_4$.

Case 3: $\textbf{x}=(a+1,s,a+1,u)$ where $0 \leq s \leq a$ and $0 \leq u \leq a+1$.
Let $\textbf{x}' = \textbf{x} - \textbf{v}_8 = (1, s-a, 1,u-a-1)$, which lies between $\textbf{0}$ and $-\textbf{v}_7$ unless $s = 0$, in which case let $\textbf{x}''= \textbf{x}' + \textbf{v}_7 = (-a+2,-1,-a-1,u)$ which lies between $\textbf{0}$ and $-\textbf{v}_2$ unless $u = a+1$.
If $s = 0$ and $u = a+1$ then $\textbf{x}'=(1,-a,1,0)$ which lies between $\textbf{0}$ and $-\textbf{v}_6$.

Case 4: $\textbf{x}=(r, a+1,a+1,u)$ where $0 \leq r \leq a$ and $0 \leq u \leq a+1$.
Let $\textbf{x}' = \textbf{x} - \textbf{v}_8 = (r-a, 1, 1,u-a-1)$, which lies between $\textbf{0}$ and $\textbf{v}_1$ unless $u \leq 1$, in which case let $\textbf{x}''= \textbf{x}' - \textbf{v}_1 = (r+1,-a,-a+-1,u-2)$ which lies between $\textbf{0}$ and $-\textbf{v}_5$ unless $r = a$.
If $r = a$ and $u \leq 1$ then let $\textbf{x}'''= \textbf{x}'' + \textbf{v}_5 = (1,0,0,u+a)$ which lies between $\textbf{0}$ and $\textbf{v}_8$.

Case 5: $\textbf{x}=(a+1,s,t,u)$ where $0 \leq s,t \leq a$ and $0 \leq u \leq a+1$.
Let $\textbf{x}' = \textbf{x} - \textbf{v}_8 = (1, s-a, t-a,u-a-1)$, which lies between $\textbf{0}$ and $-\textbf{v}_5$ unless $t = 0$, in which case let $\textbf{x}''= \textbf{x}' + \textbf{v}_5 = (-a+1,s,-1,u+1)$ which lies between $\textbf{0}$ and $\textbf{v}_7$ unless $s = a$ or $u=a+1$.
If $s = a$ and $t = 0$ then $\textbf{x}'= (1,0,-a,u-a-1)$ which lies between $\textbf{0}$ and $-\textbf{v}_4$.
If $t = 0$ and $u = a+1$ then $\textbf{x}'= (1,s-a,-a,0)$ which lies between $\textbf{0}$ and $-\textbf{v}_1$.

Case 6: $\textbf{x}=(r,a+1,t,u)$ where $0 \leq r,t \leq a$ and $0 \leq u \leq a+1$.
Let $\textbf{x}' = \textbf{x} - \textbf{v}_8 = (r-a, 1, t-a,u-a-1)$, which lies between $\textbf{0}$ and $\textbf{v}_6$ unless $u = 0$, in which case $\textbf{x}$ lies between $\textbf{0}$ and $\textbf{v}_2$ unless $r = a$.
If $r = a$ and $u = 0$ then let $\textbf{x}'' = \textbf{x} - \textbf{v}_2 = (1,0,t-a-1,a)$ which lies between $\textbf{0}$ and $-\textbf{v}_3$.

Case 7: $\textbf{x}=(r,s,a+1,u)$ where $0 \leq r,s \leq a$ and $0 \leq u \leq a+1$.
Let $\textbf{x}' = \textbf{x} - \textbf{v}_8 = (r-a, s-a, 1,u-a-1)$, which lies between $\textbf{0}$ and $\textbf{v}_3$ unless $s = 0$ or $u=0$, in which case let $\textbf{x}''= \textbf{x}' - \textbf{v}_3 = (r+1,s-1,-a,u-1)$.
If $s = 0$ and $u = 0$ then $\textbf{x}$ lies between $\textbf{0}$ and $-\textbf{v}_6$.
If $s = 0$ and $u \geq 1$ then $\textbf{x}''$ lies between $\textbf{0}$ and $-\textbf{v}_1$ unless $u=a+1$, in which case $\textbf{x}'$ lies between $\textbf{0}$ and $\textbf{v}_4$. 
If $s \geq 1$ and $u = 0$ then $\textbf{x}$ lies between $\textbf{0}$ and $\textbf{v}_2$ unless $r=a$, in which case $\textbf{x}'$ lies between $\textbf{0}$ and $-\textbf{v}_7$. \newline
This completes the cases for the orthant of $\textbf{v}_8$.

This also completes the proof of the theorem for any $k \equiv 0 \pmod 2$. 

Now we consider the eight orthants $\textbf{v}_1,...,\textbf{v}_8$ in turn for the case $k \equiv 1 \pmod 2$.

\section*{Orthant of $\textbf{v}_1$, $k \equiv 1$ (mod 2) }
First suppose that  $\textbf{x}$ lies within the orthant of $\textbf{v}_1$, but not between $\textbf{0}$ and $\textbf{v}_1$. Then the first coordinate of $\textbf{x}$ is equal to $-a$ or $-a-1$, or the third coordinate equals $-a$ or $-a-1$, or the fourth equals $a+1$. We distinguish seven cases.

Case 1: $\textbf{x}=(-r,s,-t,a+1)$ where $a \leq r, t \leq a+1$ and $0 \leq s \leq a+1$.
Let $\textbf{x}' = \textbf{x} - \textbf{v}_1 = (a-1-r, s-a-1, a-1-t,1)$, which lies between $\textbf{0}$ and $\textbf{v}_8$ unless $s \leq 1$ in which case let $\textbf{x}'' = \textbf{x}'-\textbf{v}_8=(2a-2-r, s-2,2a-1-t,-a)$ which lies between $\textbf{0}$ and $-\textbf{v}_1$.

Case 2: $\textbf{x}=(-r,s,-t,u)$ where $a \leq r, t\leq a+1$ and $0 \leq s \leq a+1$ and $0 \leq u \leq a$.
Let $\textbf{x}' = \textbf{x} - \textbf{v}_1 = (a-1-r, s-a-1, a-1-t, u-a)$, which lies between $\textbf{0}$ and $\textbf{v}_5$ unless $s = 0$ or $u \leq 1$, in which case let $\textbf{x}'' = \textbf{x}'-\textbf{v}_5=(2a-1-r,s-1,2a-t,u-2)$. If $s = 0$ and $u \leq 1$ then $\textbf{x}''$ lies between $\textbf{0}$ and $-\textbf{v}_1$, unless $t=a$, in which case let $\textbf{x}'''=\textbf{x}''+\textbf{v}_1=(a-r,a,1,u+a-2)$ which lies between $\textbf{0}$ and $\textbf{v}_4$.
If $s = 0$ and $u \geq 2$ then $\textbf{x}''$ lies between $\textbf{0}$ and $-\textbf{v}_7$. If $s \geq 1$ and $u \leq 1$ then $\textbf{x}''$ lies between $\textbf{0}$ and $-\textbf{v}_8$ unless $s=a+1$, in which case let $\textbf{x}''''=\textbf{x}''+\textbf{v}_8=(a-r,1,a-t,a+u-1)$ which lies between $\textbf{0}$ and $\textbf{v}_1$.

Case 3: $\textbf{x}=(-r,s,-t,a+1)$ where $a\leq r\leq a+1$, $0\leq s \leq a+1$ and $0 \leq t \leq a-1$.
Let $\textbf{x}' = \textbf{x} - \textbf{v}_1 = (a-1-r, s-a-1, a-1-t,1)$, which lies between $\textbf{0}$ and $-\textbf{v}_6$ unless $s = 0$, in which case let $\textbf{x}'' = \textbf{x}'+\textbf{v}_6 = (2a-1-r,-1,-t,-a+1)$  which lies between $\textbf{0}$ and $-\textbf{v}_4$.

Case 4: $\textbf{x}=(-r,s,-t,a+1)$ where $0\leq r\leq a-1$, $0\leq s \leq a+1$ and $a \leq t \leq a+1$.
Let $\textbf{x}' = \textbf{x} - \textbf{v}_1 = (a-1-r, s-a-1, a-1-t,1)$, which lies between $\textbf{0}$ and $-\textbf{v}_3$ unless $s \leq 1$, in which case let $\textbf{x}'' = \textbf{x}'+\textbf{v}_3 = (-2-r,s-2,2a-2-t,-a+1)$  which lies between $\textbf{0}$ and $-\textbf{v}_2$.

Case 5: $\textbf{x}=(-r,s,-t,a+1)$ where $0\leq r,t\leq a-1$ and $0 \leq s \leq a+1$.
Let $\textbf{x}' = \textbf{x} - \textbf{v}_1 = (a-1-r, s-a-1, a-1-t,1)$, which lies between $\textbf{0}$ and $-\textbf{v}_7$ unless $s =0$, in which case let $\textbf{x}'' = \textbf{x}'+\textbf{v}_7 = (-r-1,-1,-t-1,-a+2)$  which lies between $\textbf{0}$ and $\textbf{v}_5$.

Case 6: $\textbf{x}=(-r,s,-t,u)$ where $0\leq r\leq a-1$, $0\leq s \leq a+1$, $a \leq t \leq a+1$ and $0 \leq u \leq a$.
Let $\textbf{x}' = \textbf{x} - \textbf{v}_1 = (a-1-r, s-a-1, a-1-t,u-a)$, which lies between $\textbf{0}$ and $-\textbf{v}_4$ unless $s = 0$ or $u=0$, in which case let $\textbf{x}'' = \textbf{x}'+\textbf{v}_4 = (-r-1,s-1,2a-1-t,u-1)$. If $s=0$ and $u=0$ then let $\textbf{x}''' = \textbf{x}''+\textbf{v}_2=(a-r,a,a+1-t,a-2)$ which lies between $\textbf{0}$ and $-\textbf{v}_5$.
If $s = 0$ and $u \geq 1$ then $\textbf{x}''$ lies between $\textbf{0}$ and $-\textbf{v}_6$.
If $s \geq 1$ and $u =0$ then $\textbf{x}''$ lies between $\textbf{0}$ and $\textbf{v}_3$ unless $s=a+1$, in which case $\textbf{x}'$ lies between $\textbf{0}$ and $\textbf{v}_6$.

Case 7: $\textbf{x}=(-r,s,-t,u)$ where $a \leq r\leq a+1$, $0\leq s \leq a+1$, $0 \leq t \leq a-1$ and $0 \leq u \leq a$.
Let $\textbf{x}' = \textbf{x} - \textbf{v}_1 = (a-1-r, s-a-1, a-1-t,u-a)$, which lies between $\textbf{0}$ and $-\textbf{v}_2$ unless $t = 0$ or $u=0$, in which case let $\textbf{x}'' = \textbf{x}'+\textbf{v}_2 = (2a-r,s,-t+1,u-1)$. If $t=0$ and $u=0$ then let $\textbf{x}''' = \textbf{x}''+\textbf{v}_8=(a+1-r,s-a+1,-a+1,a)$ which lies between $\textbf{0}$ and $-\textbf{v}_3$ unless $a \leq s \leq a+1$, in which case let $\textbf{x}'''' = \textbf{x}-\textbf{v}_7=(a-r,s-a,a,a-1)$ which lies between $\textbf{0}$ and $\textbf{v}_4$.
If $t = 0$ and $u \geq 1$ then $\textbf{x}''$ lies between $\textbf{0}$ and $-\textbf{v}_5$ unless $s=a+1$ or $u=a$ in which case let $\textbf{x}^v = \textbf{x}''+\textbf{v}_5=(a-r,s-a,-a,-a+u+1)$. 
If $s = a+1$ then $\textbf{x}'$ lies between $\textbf{0}$ and $\textbf{v}_3$. If $1 \leq s \leq a$ and $u=a$ then $\textbf{x}^v$ lies between $\textbf{0}$ and $\textbf{v}_8$.
If $s = 0$ and $u=a$ then $\textbf{x}''$ lies between $\textbf{0}$ and $-\textbf{v}_7$. If $t \geq 1$ and $u=0$ then $\textbf{x}''$ lies between $\textbf{0}$ and $\textbf{v}_6$ unless $s=a+1$, in which case  $\textbf{x}'$ lies between $\textbf{0}$ and $\textbf{v}_3$.

This completes the cases for the orthant of $\textbf{v}_1$. 

\section*{Orthant of $\textbf{v}_2$, $k \equiv 1$ (mod 2) }
Now suppose that $\textbf{x}$ lies in the orthant of $\textbf{v}_2$ but not between $\textbf{0}$ and $\textbf{v}_2$. Then the third coordinate of $\textbf{x}$ is equal to $-a+1,-a$ or $-a-1$, or the fourth coordinate equals $a$ or $a+1$. We distinguish three cases.

Case 1: $\textbf{x}=(r,s,-t,u)$ where $0 \leq r,s \leq a+1$, $a-1 \leq t \leq a+1$ and $a \leq u \leq a+1$.
Let $\textbf{x}' = \textbf{x} - \textbf{v}_2 = (r-a-1, s-a-1, a-2-t, u-a+1)$, which lies between $\textbf{0}$ and $\textbf{v}_8$ unless $r \leq 1$ or $s \leq 1$, in which case let $\textbf{x}'' = \textbf{x}'-\textbf{v}_8=(r-2, s-2, 2a-2-t, u-2a)$.
If $r \leq 1$ and $s \geq 2$ then $\textbf{x}''$ lies between $\textbf{0}$ and $\textbf{v}_3$.
If $r \geq 2$ and $s \leq 1$ then $\textbf{x}''$ lies between $\textbf{0}$ and $-\textbf{v}_1$.
If $r \leq 1$ and $s \leq 1$ then $\textbf{x}''$ lies between $\textbf{0}$ and $-\textbf{v}_2$ unless $t=a-1$ or $u=a$. Let $\textbf{x}'''= \textbf{x}''+ \textbf{v}_2 = (r+a-1,s+a-1,a-t,u-a-1)$. 
If $r \leq 1$ and $s \leq 1$ and $u=a$ then $\textbf{x}'''$ lies between $\textbf{0}$ and $\textbf{v}_6$ unless $t=a-1$, in which case let $\textbf{x}''''= \textbf{x}'''- \textbf{v}_6 = (r-1,s-1,2a-1-t,u-1)$ which lies between $\textbf{0}$ and $\textbf{v}_4$ if $s=1$, and between $\textbf{0}$ and $-\textbf{v}_7$ if $r=1$.
If $r =0$ and $s =0$ then $\textbf{x}$ lies between $\textbf{0}$ and $\textbf{v}_1$.
If $r \leq 1$ and $s \leq 1$ and $t = a-1$ and $u=a+1$ then $\textbf{x}'''$ lies between $\textbf{0}$ and $-\textbf{v}_5$.

Case 2: $\textbf{x}=(r,s,-t,u)$ where $0 \leq r,s \leq a+1$, $0 \leq t \leq a-2$ and $a \leq u \leq a+1$.
Let $\textbf{x}' = \textbf{x} - \textbf{v}_2 = (r-a-1,s-a-1, a-2-t, u-a+1)$, which lies between $\textbf{0}$ and $-\textbf{v}_6$ unless $r = 0$ or $s = 0$.
Let $\textbf{x}''$ = $\textbf{x}' + \textbf{v}_6 = (r-1, s-1,-t-1,u-2a+1)$. 
If $r = 0$ and $s = 0$ then $\textbf{x}''$ lies between $\textbf{0}$ and $ \textbf{v}_5$ unless $u=a$, in which case let $\textbf{x}'''$ = $\textbf{x}'' - \textbf{v}_5 = (r+a-1, s+a-1,a-t,u-a-1)$ which lies between $\textbf{0}$ and $-\textbf{v}_8$.
If $r=0$ and $s \geq 1$ then $\textbf{x}''$ lies between $\textbf{0}$ and $\textbf{v}_7$.
If $r \geq 1$ and $s=0$ then $\textbf{x}''$ lies between $\textbf{0}$ and $- \textbf{v}_4$.

Case 3: $\textbf{x}=(r,s,-t,u)$ where $0 \leq r,s \leq a+1$, $a-1 \leq t \leq a+1$ and $0 \leq u \leq a-1$.
Let $\textbf{x}' = \textbf{x} - \textbf{v}_2 = (r-a-1, s-a-1, a-2-t, u-a+1)$, which lies between $\textbf{0}$ and $\textbf{v}_5$ unless $r = 0$ or $s = 0$ or $u=0$, in which case let $\textbf{x}'' = \textbf{x}' - \textbf{v}_5 = (r-1,s-1,2a-1-t,u-1)$.
If $r=0$ and $s=0$ and $u=0$ then $\textbf{x}''$ lies between $\textbf{0}$ and $-\textbf{v}_2$ unless $a-1 \leq t \leq a$, in which case let $\textbf{x}''' = \textbf{x}'' + \textbf{v}_2 = (a+r,a+s,a+1-t,a-2+u)$ which lies between $\textbf{0}$ and $-\textbf{v}_5$.
If $r=0$ and $s=0$ and $u \geq 1$ then $\textbf{x}''$ lies between $\textbf{0}$ and $-\textbf{v}_6$ unless $t=a-1$, in which case let $\textbf{x}'''' = \textbf{x}''+\textbf{v}_6=(r+a-1,s+a-1,a-t,u-a-1)$ which lies between $\textbf{0}$ and $-\textbf{v}_8$. 
If $r=0$ and $s \geq 1$ and $u = 0$ then $\textbf{x}''$ lies between $\textbf{0}$ and $\textbf{v}_3$ unless $s=a+1$ or $t=a-1$, in which case let $\textbf{x}^v = \textbf{x}''-\textbf{v}_3=(r+a,s-a,a-t,u+a-1)$. 
If $s=a+1$ then $\textbf{x}^v$ lies between $\textbf{0}$ and $\textbf{v}_2$ unless $t=a-1$.
If $t=a-1$ then $\textbf{x}^v$ lies  between $\textbf{0}$ and $-\textbf{v}_7$ unless $s=a+1$, in which case let $\textbf{x}^{vi} = \textbf{x}^v+\textbf{v}_5=(r,s-2a,-t-1,u+1)$ which lies between $\textbf{0}$ and $\textbf{v}_8$. 
If $r=0$ and $s \geq 1$ and $u \geq 1$ then $\textbf{x}''$ lies between $\textbf{0}$ and $\textbf{v}_4$.
If $r \geq 1$ and $s=0$ and $u=0$ then $\textbf{x}''$ lies between $\textbf{0}$ and $-\textbf{v}_1$ unless $r=a+1$ or $t=a-1$, in which case let $\textbf{x}''' = \textbf{x}'' + \textbf{v}_1 = (r-a,s+a,a-t,a+u-1)$.
If $r=a+1$ and $t \geq a$ then $\textbf{x}'''$ lies between $\textbf{0}$ and $\textbf{v}_2$.
If $r=a+1$ and $t=a-1$ then let $\textbf{x}'''' = \textbf{x}''' + \textbf{v}_5 = (r-2a,s,-t-1,u+1)$ which lies between $\textbf{0}$ and $\textbf{v}_8$.
If $1 \leq r \leq a$ and $t=a-1$ then $\textbf{x}'''$ lies between $\textbf{0}$ and $\textbf{v}_4$.
If $r \geq 1$ and $s=0$ and $u \geq 1$ then $\textbf{x}''$ lies between $\textbf{0}$ and $-\textbf{v}_7$.
If $r \geq 1$ and $s \geq 1$ and $u=0$ then $\textbf{x}''$ lies between $\textbf{0}$ and $-\textbf{v}_8$ unless $r=a+1$ or $s=a+1$, in which case let $\textbf{x}''' = \textbf{x}'' + \textbf{v}_8 = (r-a,s-a,a-1-t,u+a)$.
If $r=a+1$ and $s=a+1$ then let $\textbf{x}'''' = \textbf{x}''' - \textbf{v}_2 = (r-2a-1,s-2a-1,2a-3-t,u+1)$ which lies between $\textbf{0}$ and $-\textbf{v}_6$.
If $r=a+1$ and $1 \leq s \leq a$ then $\textbf{x}'''$ lies between $\textbf{0}$ and $-\textbf{v}_3$.
If $1 \leq r \leq a$ and $s=a+1$ then $\textbf{x}'''$ lies between $\textbf{0}$ and $\textbf{v}_1$.

This completes the cases for the orthant of $\textbf{v}_2$.

\section*{Orthant of $\textbf{v}_3$, $k \equiv 1$ (mod 2) }
Now suppose that $\textbf{x}$ lies in the orthant of $\textbf{v}_3$ but not between $\textbf{0}$ and $\textbf{v}_3$. Then the second coordinate of $\textbf{x}$ is equal to $a$ or $a+1$, or the third coordinate equals $a$ or $a+1$, or the fourth equals $-a-1$. We distinguish seven cases.

Case 1: $\textbf{x}=(-r,s,t,-a-1)$ where $0 \leq r \leq a+1$ and $a \leq s,t \leq a+1$.
Let $\textbf{x}' = \textbf{x} - \textbf{v}_3 = (a+1-r, s-a+1, t-a+1, -1)$, which lies between $\textbf{0}$ and $-\textbf{v}_8$ unless $r \leq 1$, in which case let $\textbf{x}'' = \textbf{x}'+\textbf{v}_8=(2-r, s-2a+2, t-2a+1, a)$ which lies between $\textbf{0}$ and $-\textbf{v}_3$.

Case 2: $\textbf{x}=(-r,s,t,-u)$ where $0 \leq r \leq a+1$ and $a \leq s,t \leq a+1$ and $0 \leq u \leq a$.
Let $\textbf{x}' = \textbf{x} - \textbf{v}_3 = (a+1-r, s-a+1, t-a+1, a-u)$, which lies between $\textbf{0}$ and $-\textbf{v}_5$ unless $r = 0$ or $u \leq 1$, in which case let $\textbf{x}'' = \textbf{x}'+\textbf{v}_5=(1-r,s-2a+1,t-2a,2-u)$.
If $r = 0$ and $u \leq 1$ then $\textbf{x}''$ lies between $\textbf{0}$ and $-\textbf{v}_3$ unless $t=a$, in which case let $\textbf{x}''' = \textbf{x}''+\textbf{v}_3=(-a-r, s-a,t-a-1,2-a-u)$ which lies between $\textbf{0}$ and $\textbf{v}_7$.
If $r = 0$ and $u \geq 2$ then $\textbf{x}''$ lies between $\textbf{0}$ and $-\textbf{v}_4$.
If $r \geq 1$ and $u \leq 1$ then $\textbf{x}''$ lies between $\textbf{0}$ and $\textbf{v}_8$ unless $r=a+1$, in which case let $\textbf{x}'''' = \textbf{x}''-\textbf{v}_8=(a-r,s-a,t-a,1-a-u)$ which lies between $\textbf{0}$ and $\textbf{v}_3$.

Case 3: $\textbf{x}=(-r,s,t,-a-1)$ where $0 \leq r \leq a+1$, $a \leq s \leq a+1$ and $0 \leq t \leq a-1$.
Let $\textbf{x}' = \textbf{x} - \textbf{v}_3 = (a+1-r, s-a+1,t-a+1, -1)$, which lies between $\textbf{0}$ and $\textbf{v}_6$ unless $r = 0$, in which case let $\textbf{x}'' = \textbf{x}' - \textbf{v}_6 = (1-r, s-2a+1,t,a-1)$ which lies between $\textbf{0}$ and $-\textbf{v}_7$. 

Case 4: $\textbf{x}=(-r,s,t,-a-1)$ where $0 \leq r \leq a+1$ and $0 \leq s \leq a-1$ and $a \leq t \leq a+1$.
Let $\textbf{x}' = \textbf{x} - \textbf{v}_3 = (a+1-r, s-a+1,t-a+1,-1)$, which lies between $\textbf{0}$ and $- \textbf{v}_1$ unless $r \leq 1$, in which case let $\textbf{x}'' = \textbf{x}' + \textbf{v}_1 = (2-r, s+2,t-2a+2,a-1)$ which lies between $\textbf{0}$ and $\textbf{v}_2$. 

Case 5: $\textbf{x}=(-r,s,t,-a-1)$ where $0 \leq r \leq a+1$ and $0 \leq s,t \leq a-1$.
Let $\textbf{x}' = \textbf{x} - \textbf{v}_3 = (a+1-r, s-a+1,t-a+1,-1)$, which lies between $\textbf{0}$ and $- \textbf{v}_4$ unless $r =0$, in which case let $\textbf{x}'' = \textbf{x}' + \textbf{v}_4 = (1-r, s+1,t+1,a-2)$ which lies between $\textbf{0}$ and $-\textbf{v}_5$. 

Case 6: $\textbf{x}=(-r,s,t,-u)$ where $0 \leq r \leq a+1$, $0 \leq s \leq a-1$, $a \leq t \leq a+1$ and $0 \leq u \leq a$.
Let $\textbf{x}' = \textbf{x} - \textbf{v}_3 = (a+1-r, s-a+1,t-a+1,a-u)$, which lies between $\textbf{0}$ and $- \textbf{v}_7$ unless $r = 0$ or $u=0$, in which case let $\textbf{x}'' = \textbf{x}' + \textbf{v}_7 = (1-r, s+1,t-2a+1,1-u)$.
If $r = 0$ and $u = 0$ then $\textbf{x}''$ lies between $\textbf{0}$ and $\textbf{v}_2$ unless $t=a$, in which case let $\textbf{x}''' = \textbf{x}''-\textbf{v}_2=(-a-r, s-a,t-a-1,2-a-u)$ which lies between $\textbf{0}$ and $\textbf{v}_5$.
If $r = 0$ and $u \geq 1$ then $\textbf{x}''$ lies between $\textbf{0}$ and $\textbf{v}_6$.
If $r \geq 1$ and $u =0$ then $\textbf{x}''$ lies between $\textbf{0}$ and $\textbf{v}_1$ unless $r=a+1$, in which case let $\textbf{x}''' = \textbf{x}''-\textbf{v}_1=(a-r,s-a,t-a,1-a-u)$ which lies between $\textbf{0}$ and $-\textbf{v}_2$.

Case 7: $\textbf{x}=(-r,s,t,-u)$ where $0 \leq r \leq a+1$, $a \leq s \leq a+1$, $0 \leq t \leq a-1$ and $0 \leq u \leq a$.
Let $\textbf{x}' = \textbf{x} - \textbf{v}_3 = (a+1-r, s-a+1,t-a+1,a-u)$, which lies between $\textbf{0}$ and $ \textbf{v}_2$ unless $t = 0$ or $u=0$, in which case let $\textbf{x}'' = \textbf{x}' - \textbf{v}_2 = (-r, s-2a,t-1,1-u)$.
If $t = 0$ and $u = 0$ then $\textbf{x}''$ lies between $\textbf{0}$ and $\textbf{v}_8$ unless $r \leq 1$, in which case let $\textbf{x}''' = \textbf{x}''-\textbf{v}_8=(a-1-r, s-a-1,t+a-1,-a-u)$ which lies between $\textbf{0}$ and $-\textbf{v}_1$.
If $t = 0$ and $u \geq 1$ then $\textbf{x}''$ lies between $\textbf{0}$ and $\textbf{v}_5$ unless $r=0$ or $u=a$, in which case let $\textbf{x}'''' = \textbf{x}''-\textbf{v}_5=(a-r,s-a,a+t,a-u-1)$.
If $r = 0$, $t=0$ and $u =a$ then let $\textbf{x}^v = \textbf{x}''''+\textbf{v}_8=(1-r,s-2a+1,t,2a-u)$ which lies between $\textbf{0}$ and $-\textbf{v}_3$.
If $r = 0$, $t=0$ and $1 \leq u \leq a-1$ then $\textbf{x}''''$ lies between $\textbf{0}$ and $-\textbf{v}_5$.
If $1 \leq r \leq a$, $t=0$ and $u =a$ then $\textbf{x}'''$ lies between $\textbf{0}$ and $-\textbf{v}_8$.
If $r = a+1$, $t=0$ and $u =a$ then $\textbf{x}^v$ lies between $\textbf{0}$ and $-\textbf{v}_6$.
If $t \geq 1$ and $u =0$ then $\textbf{x}''$ lies between $\textbf{0}$ and $-\textbf{v}_6$.

This completes the cases for the orthant of $\textbf{v}_3$.

\section*{Orthant of $\textbf{v}_4$, $k \equiv 1$ (mod 2) }
Now suppose $\textbf{x}$ lies in the orthant of $\textbf{v}_4$ but not between $\textbf{0}$ and $\textbf{v}_4$. Then the first coordinate of $\textbf{x}$ is equal to $-a-1$ or the second coordinate is equal to $a+1$, or the third equals $a+1$ or the fourth equals $a$ or $a+1$. We distinguish fifteen cases.

Case 1: $\textbf{x}=(-a-1,a+1,a+1,u)$ where $a \leq u \leq a+1$.
Let $\textbf{x}' = \textbf{x} - \textbf{v}_4 = (-1, 1, 1,u-a+1)$, which lies between $\textbf{0}$ and $\textbf{v}_4$.

Case 2: $\textbf{x}=(-a-1,a+1,a+1,u)$ where $0 \leq u \leq a-1$.
Let $\textbf{x}' = \textbf{x} - \textbf{v}_4 = (-1, 1, 1,u-a+1)$, which lies between $\textbf{0}$ and $\textbf{v}_3$.

Case 3: $\textbf{x}=(-a-1,a+1,t,u)$ where $0 \leq t \leq a$ and $a \leq u \leq a+1$.
Let $\textbf{x}' = \textbf{x} - \textbf{v}_4 = (-1, 1, t-a,u-a+1)$, which lies between $\textbf{0}$ and $\textbf{v}_7$.

Case 4: $\textbf{x}=(-a-1,s,a+1,u)$ where $0 \leq s \leq a$ and $a \leq u \leq a+1$.
Let $\textbf{x}' = \textbf{x} - \textbf{v}_4 = (-1, s-a, 1,u-a+1)$, which lies between $\textbf{0}$ and $-\textbf{v}_2$.

Case 5: $\textbf{x}=(-r,a+1,a+1,u)$ where $0 \leq r \leq a$ and $a \leq u \leq a+1$.
Let $\textbf{x}' = \textbf{x} - \textbf{v}_4 = (a-r, 1, 1,u-a+1)$, which lies between $\textbf{0}$ and $-\textbf{v}_5$.

Case 6: $\textbf{x}=(-r,s,a+1,u)$ where $0 \leq r,s \leq a$ and $a \leq u \leq a+1$.
Let $\textbf{x}' = \textbf{x} - \textbf{v}_4 = (a-r, s-a, 1,u-a+1)$, which lies between $\textbf{0}$ and $-\textbf{v}_7$.

Case 7: $\textbf{x}=(-r,a+1,t,u)$ where $0 \leq r,t \leq a$ and $a \leq u \leq a+1$.
Let $\textbf{x}' = \textbf{x} - \textbf{v}_4 = (a-r, 1, t-a, u-a+1)$, which lies between $\textbf{0}$ and $\textbf{v}_2$ unless $t \leq 1$, in which case let $\textbf{x}''= \textbf{x}'- \textbf{v}_2 = (-r-1,-a,t-2,u)$ and $\textbf{x}'''= \textbf{x}''- \textbf{v}_8 = (a-r-2,-1,a+t-2,u-a-1)$.
$\textbf{x}'''$ lies between $\textbf{0}$ and $-\textbf{v}_1$ unless $a-1 \leq r \leq a$, in which case let $\textbf{x}''''= \textbf{x}'''+ \textbf{v}_2 = (2a-r-1,a,t,u-2)$, which lies between $\textbf{0}$ and $-\textbf{v}_5$ unless $u = a+1$.
If $a-1 \leq r \leq a$,  $t \leq 1$ and $u = a+1$, then let $\textbf{x}^v = \textbf{x}''''+\textbf{v}_5=(a-r-1,0,t-a-1,u-a)$ and $\textbf{x}^{vi} = \textbf{x}^v-\textbf{v}_8=(2a-r-2,a-1,t-1,u-2a-1)$ which lies between $\textbf{0}$ and $\textbf{v}_6$.

Case 8: $\textbf{x}=(-a-1,s,t,u)$ where $0 \leq s,t \leq a$ and $a \leq u \leq a+1$.
Let $\textbf{x}' = \textbf{x} - \textbf{v}_4 = (-1, s-a, t-a, u-a+1)$, which lies between $\textbf{0}$ and $\textbf{v}_8$ unless $s=0$, in which case let $\textbf{x}''= \textbf{x}'- \textbf{v}_8 = (a-2,s-1,t,u-2a)$ which lies between $\textbf{0}$ and $-\textbf{v}_1$ unless $t=a$.
If $s=0$ and $t=a$ then let $\textbf{x}''' = \textbf{x}''+ \textbf{v}_1 = (-1,s+a,t-a+1,u-a)$, which lies between $\textbf{0}$ and $\textbf{v}_4$.

Case 9: $\textbf{x}=(-r,a+1,a+1,u)$ where $0 \leq r \leq a$ and $0 \leq u \leq a-1$.
Let $\textbf{x}' = \textbf{x} - \textbf{v}_4 = (a-r,1,1, u-a+1)$, which lies between $\textbf{0}$ and $-\textbf{v}_8$ unless $r=0$, in which case let $\textbf{x}''= \textbf{x}'+ \textbf{v}_8 = (1-r,2-a,1-a,u+2)$ which lies between $\textbf{0}$ and $\textbf{v}_8$.

Case 10: $\textbf{x}=(-a-1,s,a+1,u)$ where $0 \leq s \leq a$ and $0 \leq u \leq a-1$.
Let $\textbf{x}' = \textbf{x} - \textbf{v}_4 = (-1,s-a,1,u-a+1)$, which lies between $\textbf{0}$ and $-\textbf{v}_2$.

Case 11: $\textbf{x}=(-a-1,a+1,t,u)$ where $0 \leq t \leq a$ and $0 \leq u \leq a-1$.
Let $\textbf{x}' = \textbf{x} - \textbf{v}_4 = (-1,1,t-a,u-a+1)$, which lies between $\textbf{0}$ and $\textbf{v}_7$.

Case 12: $\textbf{x}=(-r,s,t,u)$ where $0 \leq r,s,t \leq a$ and $a \leq u \leq a+1$.
Let $\textbf{x}' = \textbf{x} - \textbf{v}_4 = (a-r,s-a,t-a,u-a+1)$, which lies between $\textbf{0}$ and $-\textbf{v}_3$ unless $s=0$ or $t=0$.
If $s = 0$ and $t=0$ then let $\textbf{x}'' = \textbf{x}'+\textbf{v}_3=(-1-r,s-1,t-1,u-2a+1)$ which lies between $\textbf{0}$ and $\textbf{v}_5$ unless $r=a$ or $u=a$ in which case let $\textbf{x}''' = \textbf{x}''-\textbf{v}_5=(a-1-r,s+a-1,t+a,u-a-1)$.
If $r = a$ and $u =a$, then $\textbf{x}$ lies between $\textbf{0}$ and $-\textbf{v}_6$.
If $r = a$ and $u =a+1$, then $\textbf{x}'''$ lies between $\textbf{0}$ and $\textbf{v}_4$.
If $r \leq a-1$ and $u =a$, then $\textbf{x}'''$ lies between $\textbf{0}$ and $-\textbf{v}_8$.
If $s = 0$ and $1 \leq t \leq a$ then let $\textbf{x}''$ lies between $\textbf{0}$ and $-\textbf{v}_2$ unless $t=a$, in which case let $\textbf{x}''' = \textbf{x}''+\textbf{v}_2=(a-r,a+s,t-a+1,u-a)$ which lies between $\textbf{0}$ and $-\textbf{v}_5$.
If $1 \leq s \leq a$ and $t=0$ then $\textbf{x}''$ lies between $\textbf{0}$ and $\textbf{v}_7$ unless $r=a$, in which case let $\textbf{x}''' = \textbf{x}''-\textbf{v}_7=(a-1-r,s-a-1,t+a-1,u-a)$ which lies between $\textbf{0}$ and $-\textbf{v}_6$.

Case 13: $\textbf{x}=(-r,s,a+1,u)$ where $0 \leq r,s \leq a$ and $0 \leq u \leq a-1$.
Let $\textbf{x}' = \textbf{x} - \textbf{v}_4 = (a-r, s-a, 1,u-a+1)$, which lies between $\textbf{0}$ and $-\textbf{v}_1$ unless $r = 0$, in which case let $\textbf{x}'' = \textbf{x}'+\textbf{v}_1=(1-r, s+1,-a+2,u+1)$ which lies between $\textbf{0}$ and $\textbf{v}_2$ unless $u = a-1$.
If $r = 0$ and $u = a-1$ then let $\textbf{x}''' = \textbf{x}''-\textbf{v}_2=(-a-r,s-a,0,u-a+2)$ which lies between $\textbf{0}$ and $\textbf{v}_5$.

Case 14: $\textbf{x}=(-r,a+1, t,u)$ where $0 \leq r,t \leq a$ and $0 \leq u \leq a-1$.
Let $\textbf{x}' = \textbf{x} - \textbf{v}_4 = (a-r,1,t-a,u-a+1)$, which lies between $\textbf{0}$ and $\textbf{v}_6$ unless $t = 0$, in which case let $\textbf{x}'' = \textbf{x}'-\textbf{v}_6=(-r,1-a,t-1,u+1)$ which lies between $\textbf{0}$ and $\textbf{v}_8$ unless $r = a$. 
If $r = a$ and $t=0$ then let $\textbf{x}''' = \textbf{x}''-\textbf{v}_8=(a-1-r,0,t+a-1,u-a)$ which lies between $\textbf{0}$ and $\textbf{v}_3$.

Case 15: $\textbf{x}=(-a-1,s,t,u)$ where $0 \leq s,t \leq a$ and $0 \leq u \leq a-1$.
Let $\textbf{x}' = \textbf{x} - \textbf{v}_4 = (-1,s-a,t-a,u-a+1)$, which lies between $\textbf{0}$ and $\textbf{v}_5$ unless $t=0$ or $u=0$ in which case 
let $\textbf{x}'' = \textbf{x}'-\textbf{v}_5=(a-1,s,t+1,u-1)$.
If $t = 0$ and $u =0$, then $\textbf{x}''$ lies between $\textbf{0}$ and $-\textbf{v}_8$ unless $s=a$, in which case let $\textbf{x}''' = \textbf{x}''+\textbf{v}_8=(0,s-a+1,t-a+1,u+a)$ which lies between $\textbf{0}$ and $\textbf{v}_1$.
If $t = 0$ and $1 \leq u \leq  a-1$, then $\textbf{x}''$ lies between $\textbf{0}$ and $-\textbf{v}_5$.
If $1 \leq t \leq a-1$ and $u =0$, then $\textbf{x}''$ lies between $\textbf{0}$ and $-\textbf{v}_8$ unless $s=a$, in which case $\textbf{x}'''$  lies between $\textbf{0}$ and $\textbf{v}_1$.
If $t = a$ and $u = 0$ then $\textbf{x}'''$ lies between $\textbf{0}$ and $-\textbf{v}_6$ unless $s=a$, in which case let $\textbf{x}'''' = \textbf{x}'''-\textbf{v}_4=(a,s-2a+1,t-2a+1,u+1)$ which lies between $\textbf{0}$ and $-\textbf{v}_3$.

This completes the cases for the orthant of $\textbf{v}_4$.

\section*{Orthant of $\textbf{v}_5$, $k \equiv 1$ (mod 2) }
Now suppose $\textbf{x}$ lies in the orthant of $\textbf{v}_5$ but not between $\textbf{0}$ and $\textbf{v}_5$. Then the first coordinate of $\textbf{x}$ is equal to $-a-1$ or the second coordinate is equal to $-a-1$, or the fourth equals $-a+1$, $-a$ or $-a-1$. We distinguish seven cases.

Case 1: $\textbf{x}=(-a-1,-a-1,-t,-u)$ where $0 \leq t \leq a+1$ and $a-1 \leq u \leq a+1$.
Let $\textbf{x}' = \textbf{x} - \textbf{v}_5 = (-1, -1, a+1-t,a-2-u)$, which lies between $\textbf{0}$ and $-\textbf{v}_2$ unless $t \leq 2$, in which case let $\textbf{x}''= \textbf{x}'+ \textbf{v}_2 = (a,a,3-t,2a-3-u)$ which lies between $\textbf{0}$ and $-\textbf{v}_5$.

Case 2: $\textbf{x}=(-a-1,-a-1,-t,-u)$ where $0 \leq t \leq a+1$ and $0 \leq u \leq a-2$.
Let $\textbf{x}' = \textbf{x} - \textbf{v}_5 = (-1, -1, a+1-t,a-2-u)$, which lies between $\textbf{0}$ and $-\textbf{v}_6$ unless $t \leq 1$, in which case let $\textbf{x}''= \textbf{x}'+ \textbf{v}_6 = (a-1,a-1,2-t,-2-u)$ which lies between $\textbf{0}$ and $-\textbf{v}_8$.

Case 3: $\textbf{x}=(-a-1,-s,-t,-u)$ where $0 \leq s \leq a$, $0 \leq t \leq a+1$ and $a-1 \leq u \leq a+1$.
Let $\textbf{x}' = \textbf{x} - \textbf{v}_5 = (-1, a-s, a+1-t,a-2-u)$, which lies between $\textbf{0}$ and $\textbf{v}_3$ unless $s = 0$ or $t \leq 1$, in which case let $\textbf{x}''= \textbf{x}'- \textbf{v}_3 = (a,1-s,2-t,2a-2-u)$.
If $s = 0$ and $t \leq 1$ then $\textbf{x}''$ lies between $\textbf{0}$ and $-\textbf{v}_5$ unless $u=a-1$, in which case let $\textbf{x}''' = \textbf{x}''+\textbf{v}_5=(0,1-a-s,1-a-t,a-u)$ which lies between $\textbf{0}$ and $\textbf{v}_8$.
If $s = 0$ and $t \geq 2$ then $\textbf{x}''$ lies between $\textbf{0}$ and $\textbf{v}_2$ unless $t=a+1$, in which case let $\textbf{x}''' = \textbf{x}''-\textbf{v}_2=(-1,-a-s,a-t,a-1-u)$ which lies between $\textbf{0}$ and $\textbf{v}_8$.
If $1 \leq s \leq a$ and $t \leq 1$ then $\textbf{x}''$ lies between $\textbf{0}$ and $-\textbf{v}_7$.

Case 4: $\textbf{x}=(-r,-a-1,-t,-u)$ where $0 \leq r \leq a$, $0 \leq t \leq a+1$ and $a-1 \leq u \leq a+1$.
Let $\textbf{x}' = \textbf{x} - \textbf{v}_5 = (a-r, -1, a+1-t, a-2-u)$, which lies between $\textbf{0}$ and $-\textbf{v}_1$ unless $r = 0$ or $t \leq 1$, in which case let $\textbf{x}''= \textbf{x}'+ \textbf{v}_1 = (1-r,a,2-t,2a-2-u)$. 
If $r = 0$ and $t \leq 1$ then $\textbf{x}''$ lies between $\textbf{0}$ and $-\textbf{v}_5$ unless $u=a-1$, in which case let $\textbf{x}''' = \textbf{x}''+\textbf{v}_5=(1-a-r,0,1-a-t,a-u)$ which lies between $\textbf{0}$ and $\textbf{v}_8$.
If $r = 0$ and $2 \leq t \leq a+1$ then $\textbf{x}''$ lies between $\textbf{0}$ and $\textbf{v}_2$ unless $t=a+1$, in which case let $\textbf{x}''' = \textbf{x}''-\textbf{v}_2=(-a-r,-1,a-t,a-1-u)$ which lies between $\textbf{0}$ and $\textbf{v}_5$.
If $r \geq 1$ and $t \leq 1$ then $\textbf{x}''$ lies between $\textbf{0}$ and $\textbf{v}_4$.

Case 5: $\textbf{x}=(-r,-s,-t,-u)$ where $0 \leq r,s \leq a$, $0 \leq t \leq a+1$ and $a-1 \leq u \leq a+1$.
Let $\textbf{x}' = \textbf{x} - \textbf{v}_5 = (a-r, a-s, a+1-t,a-2-u)$, which lies between $\textbf{0}$ and $-\textbf{v}_8$ unless $r=0$ or $s=0$ or $t=0$,  
in which case let $\textbf{x}'' = \textbf{x}'+\textbf{v}_8=(1-r,1-s,1-t,2a-1-u)$.
If $r = 0$, $s=0$ and $t=0$ then $\textbf{x}$ lies between $\textbf{0}$ and $-\textbf{v}_8$.
If $r = 0, s=0$ and $1 \leq t \leq a+1$ then $\textbf{x}''$ lies between $\textbf{0}$ and $\textbf{v}_2$ unless $a \leq t \leq a+1$ or $u=a-1$,  
in which case let $\textbf{x}''' = \textbf{x}''-\textbf{v}_2=(-a-r,-a-s,a-1-t,a-u)$.
If $a \leq t \leq a+1$ and $a \leq u \leq a+1$ then $\textbf{x}'''$ lies between $\textbf{0}$ and $\textbf{v}_5$.
If $1 \leq t \leq a-1$ and $u = a-1$, then $\textbf{x}'''$ lies between $\textbf{0}$ and $-\textbf{v}_6$.
If $a \leq t \leq a+1$ and $u =a-1$, then let $\textbf{x}''''= \textbf{x}-\textbf{v}_7=(a-r,-a-s,a-t,a-1-u)$ which lies between $\textbf{0}$ and $-\textbf{v}_4$.

Case 6: $\textbf{x}=(-r,-a-1,-t,-u)$ where $0 \leq r \leq a$, $0 \leq t \leq a+1$ and $0 \leq u \leq a-2$.
Let $\textbf{x}' = \textbf{x} - \textbf{v}_5 = (a-r, -1, a+1-t,a-2-u)$, which lies between $\textbf{0}$ and $-\textbf{v}_7$ unless $t = 0$, in which case let $\textbf{x}''= \textbf{x}+\textbf{v}_7=(-r,a-1,1-t,-1-u)$ which lies between $\textbf{0}$ and $\textbf{v}_3$

Case 7: $\textbf{x}=(-a-1,-s,-t,-u)$ where $0 \leq s \leq a$, $0 \leq t \leq a+1$ and $0 \leq u \leq a-2$.
Let $\textbf{x}' = \textbf{x} - \textbf{v}_5 = (-1,a-s, a+1-t,a-2-u)$, which lies between $\textbf{0}$ and $\textbf{v}_4$ unless $t = 0$, in which case let $\textbf{x}''= \textbf{x}-\textbf{v}_4=(a-1,-s,1-t,-1-u)$ which lies between $\textbf{0}$ and $-\textbf{v}_1$.

This completes the cases for the orthant of $\textbf{v}_5$.

\section*{Orthant of $\textbf{v}_6$, $k \equiv 1$ (mod 2) }
Now suppose $\textbf{x}$ lies in the orthant of $\textbf{v}_6$ but not between $\textbf{0}$ and $\textbf{v}_6$. Then the first coordinate of $\textbf{x}$ is equal to $a+1$ or the second coordinate is equal to $a+1$, or the third equals $-a$ or $-a-1$ or the fourth equals $-a-1$. We distinguish fifteen cases.

Case 1: $\textbf{x}=(a+1,a+1,-t,-a-1)$ where $a \leq t \leq a+1$.
Let $\textbf{x}' = \textbf{x} - \textbf{v}_6 = (1, 1, a-1-t,-1)$, which lies between $\textbf{0}$ and $\textbf{v}_6$.

Case 2: $\textbf{x}=(a+1,a+1,-t,-u)$ where $a \leq t \leq a+1$ and $0 \leq u \leq a$.
Let $\textbf{x}' = \textbf{x} - \textbf{v}_6 = (1, 1, a-1-t,a-u)$, which lies between $\textbf{0}$ and $\textbf{v}_2$ unless $u = 0$, in which case let $\textbf{x}''= \textbf{x}'- \textbf{v}_2 = (-a,-a,2a-3-t,1-u)$ which lies between $\textbf{0}$ and $-\textbf{v}_6$.

Case 3: $\textbf{x}=(a+1,a+1,-t,-a-1)$ where $0 \leq t \leq a-1$.
Let $\textbf{x}' = \textbf{x} - \textbf{v}_6 = (1, 1, a-1-t,-1)$, which lies between $\textbf{0}$ and $-\textbf{v}_8$.

Case 4: $\textbf{x}=(a+1,s,-t,-a-1)$ where $0 \leq s \leq a$ and $a \leq t \leq a+1$.
Let $\textbf{x}' = \textbf{x} - \textbf{v}_6 = (1, s-a, a-1-t,-1)$, which lies between $\textbf{0}$ and $-\textbf{v}_4$.

Case 5: $\textbf{x}=(r,a+1,-t,-a-1)$ where $0 \leq r \leq a$ and $a \leq t \leq a+1$.
Let $\textbf{x}' = \textbf{x} - \textbf{v}_6 = (r-a,1, a-1-t,-1)$, which lies between $\textbf{0}$ and $\textbf{v}_7$.

Case 6: $\textbf{x}=(r,s,-t,-a-1)$ where $0 \leq r,s \leq a$ and $a \leq t \leq a+1$.
Let $\textbf{x}' = \textbf{x} - \textbf{v}_6 = (r-a, s-a, a-1-t,-1)$, which lies between $\textbf{0}$ and $\textbf{v}_5$.

Case 7: $\textbf{x}=(r,a+1,-t,-a-1)$ where $0 \leq r \leq a$ and $0 \leq t \leq a-1$.
Let $\textbf{x}' = \textbf{x} - \textbf{v}_6 = (r-a, 1, a-1-t,-1)$, which lies between $\textbf{0}$ and $\textbf{v}_3$.

Case 8: $\textbf{x}=(a+1,s,-t,-a-1)$ where $0 \leq s \leq a$ and $0 \leq t \leq a-1$.
Let $\textbf{x}' = \textbf{x} - \textbf{v}_6 = (1, s-a, a-1-t,-1)$, which lies between $\textbf{0}$ and $-\textbf{v}_1$.

Case 9: $\textbf{x}=(r,a+1,-t,-u)$ where $0 \leq r,u \leq a$ and $a \leq t \leq a+1$.
Let $\textbf{x}' = \textbf{x} - \textbf{v}_6 = (r-a,1,a-1-t,a-u)$, which lies between $\textbf{0}$ and $\textbf{v}_1$ unless $r = 0$, in which case let $\textbf{x}''= \textbf{x}'- \textbf{v}_1 = (r-1,-a,2a-2-t,-u)$ which lies between $\textbf{0}$ and $-\textbf{v}_2$ unless $u = a$ in which case $\textbf{x}'$ lies between $\textbf{0}$ and $\textbf{v}_7$.

Case 10: $\textbf{x}=(a+1,s,-t,-u)$ where $0 \leq s,u \leq a$ and $a \leq t \leq a+1$.
Let $\textbf{x}' = \textbf{x} - \textbf{v}_6 = (1,s-a,a-1-t,a-u)$, which lies between $\textbf{0}$ and $-\textbf{v}_3$ unless $s = 0$, in which case let $\textbf{x}''= \textbf{x}'+ \textbf{v}_3 = (-a,s-1,2a-2-t,-u)$ which lies between $\textbf{0}$ and $-\textbf{v}_2$ unless $u = a$.
If $s = 0$ and $u =a $ then let $\textbf{x}''' = \textbf{x}''+\textbf{v}_2=(1,s+a,a-t,a-1-u)$ which lies between $\textbf{0}$ and $\textbf{v}_2$.

Case 11: $\textbf{x}=(a+1,a+1,-t,-u)$ where $0 \leq t \leq a-1$ and $0 \leq u \leq a$.
Let $\textbf{x}' = \textbf{x} - \textbf{v}_6 = (1,1, a-1-t, a-u)$, which lies between $\textbf{0}$ and $-\textbf{v}_5$ unless $a-1 \leq u \leq a$, in which case let $\textbf{x}''= \textbf{x}'+ \textbf{v}_5 = (1-a,1-a,-2-t,2-u)$ which lies between $\textbf{0}$ and $\textbf{v}_5$.

Case 12: $\textbf{x}=(r,s, -t,-a-1)$ where $0 \leq r,s \leq a$ and $0 \leq t \leq a-1$.
Let $\textbf{x}' = \textbf{x} - \textbf{v}_6 = (r-a,s-a,a-1-t,-1)$, which lies between $\textbf{0}$ and $-\textbf{v}_2$ unless $t=0$, in which case let $\textbf{x}''= \textbf{x}'+ \textbf{v}_2 = (r+1,s+1,1-t,a-2)$ and $\textbf{x}'''= \textbf{x}''+ \textbf{v}_5 = (r-a+1,s-a+1,-a-t,0)$. Then  $\textbf{x}'''$ lies between $\textbf{0}$ and $\textbf{v}_8$.

Case 13: $\textbf{x}=(r,s,-t,-u)$ where $0 \leq r,s,u \leq a$ and $a \leq t \leq a+1$.
Let $\textbf{x}' = \textbf{x} - \textbf{v}_6 = (r-a, s-a, a-1-t, a-u)$, which lies between $\textbf{0}$ and $\textbf{v}_8$ unless $r = 0$ or $s = 0$, in which case let $\textbf{x}''= \textbf{x}'- \textbf{v}_8 = (r-1,s-1,2a-1-t,-1-u)$. 
If $r = 0$ and $s = 0$ then $\textbf{x}''$ lies between $\textbf{0}$ and $-\textbf{v}_2$ unless $t=a$ or $a-1 \leq u \leq a$, in which case let $\textbf{x}''' = \textbf{x}''+\textbf{v}_2=(a,a,a+1-t,a-2-u)$.
If $t =a$ and $0 \leq u \leq a-2$ then $\textbf{x}'''$ lies between $\textbf{0}$ and $-\textbf{v}_5$.
If $a \leq t \leq a+1$ and $a-1 \leq u \leq a$ then let $\textbf{x}'''' = \textbf{x}-\textbf{v}_4=(a,-a,t-a,u-a+1)$ which lies between $\textbf{0}$ and $-\textbf{v}_7$.
If $r=0$ and $1 \leq s \leq a$ then let $\textbf{x}'' = \textbf{x}-\textbf{v}_7=(r+a,s-a,a-t,a-1-u)$ which lies between $\textbf{0}$ and $-\textbf{v}_3$ unless $u=a$, in which case $\textbf{x}''$ lies between $\textbf{0}$ and $-\textbf{v}_4$.
If $1 \leq r \leq a$ and $s=0$ then let $\textbf{x}'' = \textbf{x}+\textbf{v}_4=(r-a,s+a,a-t,a-1-u)$ which lies between $\textbf{0}$ and $\textbf{v}_1$ unless $u=a$, in which case $\textbf{x}''$ lies between $\textbf{0}$ and $\textbf{v}_7$.

Case 14: $\textbf{x}=(r,a+1,-t,-u)$ where $0 \leq r,u \leq a$ and $0 \leq t \leq a-1$.
Let $\textbf{x}' = \textbf{x} - \textbf{v}_6 = (r-a, 1, a-1-t,a-u)$, which lies between $\textbf{0}$ and $\textbf{v}_4$ unless $u=0$, in which case let $\textbf{x}'' = \textbf{x}' - \textbf{v}_4 = (r,1-a,-1-t,1-u)$ which lies between $\textbf{0}$ and $-\textbf{v}_3$ unless $t=a-1$.
If $t = a-1$ and $u=0$ then let $\textbf{x}''' = \textbf{x}''+\textbf{v}_3=(r-a-1,0,a-2-t,1-a-u)$ which lies between $\textbf{0}$ and $\textbf{v}_7$ unless $r = 0$, in which case $\textbf{x}$ lies between $\textbf{0}$ and $\textbf{v}_1$.

Case 15: $\textbf{x}=(a+1,s,-t,-u)$ where $0 \leq s,u \leq a$ and $0 \leq t \leq a-1$.
Let $\textbf{x}' = \textbf{x} - \textbf{v}_6 = (1, s-a, a-1-t,a-u)$, which lies between $\textbf{0}$ and $-\textbf{v}_7$ unless $u=0$, in which case let $\textbf{x}'' = \textbf{x}' + \textbf{v}_7 = (1-a,s,-1-t,1-u)$ which lies between $\textbf{0}$ and $\textbf{v}_1$ unless $t=a-1$.
If $t = a-1$ and $u=0$ then let $\textbf{x}''' = \textbf{x}''-\textbf{v}_1=(0,s-a-1,a-2-t,1-a-u)$ which lies between $\textbf{0}$ and $-\textbf{v}_4$ unless $s = 0$, in which case $\textbf{x}$ lies between $\textbf{0}$ and $-\textbf{v}_3$.

This completes the cases for the orthant of $\textbf{v}_6$.

\section*{Orthant of $\textbf{v}_7$, $k \equiv 1$ (mod 2) }
Now suppose $\textbf{x}$ lies in the orthant of $\textbf{v}_7$ but not between $\textbf{0}$ and $\textbf{v}_7$. Then the first coordinate of $\textbf{x}$ is equal to $-a-1$ or the seond is equal to $a+1$ or the third equals $-a-1$, or the fourth equals $-a$ or $-a-1$. We distinguish fifteen cases.

Case 1: $\textbf{x}=(-a-1,a+1,-a-1,-u)$ where $a \leq u \leq a+1$.
Let $\textbf{x}' = \textbf{x} - \textbf{v}_7 = (-1,1, -1,a-1-u)$, which lies between $\textbf{0}$ and $\textbf{v}_7$.

Case 2: $\textbf{x}=(-a-1,a+1,-a-1,-u)$ where $0 \leq u \leq a-1$.
Let $\textbf{x}' = \textbf{x} - \textbf{v}_7 = (-1,1, -1,a-1-u)$, which lies between $\textbf{0}$ and $\textbf{v}_1$.

Case 3: $\textbf{x}=(-a-1,a+1,-t,-u)$ where $0 \leq t \leq a$ and $a \leq u \leq a+1$.
Let $\textbf{x}' = \textbf{x} - \textbf{v}_7 = (-1,1, a-t,a-1-u)$, which lies between $\textbf{0}$ and $\textbf{v}_3$ unless $t=0$, in which case let $\textbf{x}'' = \textbf{x}' - \textbf{v}_3 = (a,2-a, 1-t,2a-1-u)$, which lies between $\textbf{0}$ and $-\textbf{v}_7$.

Case 4: $\textbf{x}=(-a-1,s,-a-1,-u)$ where $0 \leq s \leq a$ and $a \leq u \leq a+1$.
Let $\textbf{x}' = \textbf{x} - \textbf{v}_7 = (-1,s-a, -1,a-1-u)$, which lies between $\textbf{0}$ and $\textbf{v}_5$.

Case 5: $\textbf{x}=(-r,a+1,-a-1,-u)$ where $0 \leq r \leq a$ and $a \leq u \leq a+1$.
Let $\textbf{x}' = \textbf{x} - \textbf{v}_7 = (a-r,1, -1,a-1-u)$, which lies between $\textbf{0}$ and $\textbf{v}_6$.

Case 6: $\textbf{x}=(-r,s,-a-1,-u)$ where $0 \leq r,s \leq a$ and $a \leq u \leq a+1$.
Let $\textbf{x}' = \textbf{x} - \textbf{v}_7 = (a-r,s-a, -1,a-1-u)$, which lies between $\textbf{0}$ and $-\textbf{v}_4$.

Case 7: $\textbf{x}=(-r,a+1,-t,-u)$ where $0 \leq r,t \leq a$ and $a \leq u \leq a+1$.
Let $\textbf{x}' = \textbf{x} - \textbf{v}_7 = (a-r, 1, a-t,a-1-u)$ which lies between $\textbf{0}$ and $-\textbf{v}_8$ unless $r=0$, in which case let $\textbf{x}'' = \textbf{x}' + \textbf{v}_8 = (1-r, 2-a, -t,2a-u)$ lies between $\textbf{0}$ and $-\textbf{v}_3$ unless $t=a$.
If $r=0$ and $t=a$ then let $\textbf{x}'''  = \textbf{x}'' + \textbf{v}_3 = (-a-r, 1, a-1-t,a-u)$ which lies between $\textbf{0}$ and $\textbf{v}_7$.

Case 8: $\textbf{x}=(-a-1,s,-t,-u)$ where $0 \leq s,t \leq a$ and $a \leq u \leq a+1$.
Let $\textbf{x}' = \textbf{x} - \textbf{v}_7 = (-1,s-a, a-t,a-1-u)$ lies between $\textbf{0}$ and $-\textbf{v}_2$ unless $t \leq 1$, in which case let $\textbf{x}'' = \textbf{x}' + \textbf{v}_2 = (a, s+1, 2-t,2a-2-u)$ lies between $\textbf{0}$ and $-\textbf{v}_5$ unless $s=a$.
If $s=a$ and $t \leq 1$ then let $\textbf{x}'''  = \textbf{x}'' + \textbf{v}_5 = (0, s-a+1, -a+1-t,a-u)$ which lies between $\textbf{0}$ and $\textbf{v}_7$.

Case 9: $\textbf{x}=(-r,a+1,-a-1,-u)$ where $0 \leq r \leq a$ and $0 \leq u \leq a-1$.
Let $\textbf{x}' = \textbf{x} - \textbf{v}_7 = (a-r,1, -1,a-1-u)$, which lies between $\textbf{0}$ and $\textbf{v}_2$.

Case 10: $\textbf{x}=(-a-1,s,-a-1,-u)$ where $0 \leq s \leq a$ and $0 \leq u \leq a-1$.
Let $\textbf{x}' = \textbf{x} - \textbf{v}_7 = (-1,s-a, -1,a-1-u)$ lies between $\textbf{0}$ and $\textbf{v}_8$ unless $s =0$, in which case let $\textbf{x}'' = \textbf{x}' - \textbf{v}_8 = (a-2, s-1, a-1,-2-u)$ lies between $\textbf{0}$ and $-\textbf{v}_1$ unless $u=a-1$.
If $s=0$ and $u = a-1$ then let $\textbf{x}'''  = \textbf{x}'' + \textbf{v}_1 = (-1, s+a, 0,a-2-u)$ which lies between $\textbf{0}$ and $\textbf{v}_7$.

Case 11: $\textbf{x}=(-a-1,a+1,-t,-u)$ where $0 \leq t \leq a$ and $0 \leq u \leq a-1$.
Let $\textbf{x}' = \textbf{x} - \textbf{v}_7 = (-1,1, a-t,a-1-u)$, which lies between $\textbf{0}$ and $\textbf{v}_4$.

Case 12: $\textbf{x}=(-r,s,-t,-u)$ where $0 \leq r,s,t \leq a$ and $a \leq u \leq a+1$.
Let $\textbf{x}' = \textbf{x} - \textbf{v}_7 = (a-r,s-a, a-t,a-1-u)$ lies between $\textbf{0}$ and $-\textbf{v}_1$ unless $r =0$ or $t=0$, in which case let $\textbf{x}'' = \textbf{x}' + \textbf{v}_1 = (1-r, s+1, 1-t,2a-1-u)$.
If $r=0$ and $t = 0$ then $\textbf{x}''$ lies between $\textbf{0}$ and $-\textbf{v}_5$ unless $s=a$ or $u=a$, in which case let $\textbf{x}''' = \textbf{x}'' + \textbf{v}_5 = (1-a-r, s-a+1, -a-t,a+1-u)$.
If $u=a$ then $\textbf{x}$ lies between $\textbf{0}$ and $\textbf{v}_6$.
If $s=a$ and $u=a+1$ then $\textbf{x}'''$ lies between $\textbf{0}$ and $\textbf{v}_7$.
If $r=0$ and $1 \leq t \leq a$ then $\textbf{x}''$ lies between $\textbf{0}$ and $\textbf{v}_2$ unless $t=a$, in which case $\textbf{x}'$ lies between $\textbf{0}$ and $-\textbf{v}_4$.
If $1 \leq r \leq a$ and $t = 0$ then $\textbf{x}''$ lies between $\textbf{0}$ and $\textbf{v}_4$ unless $s=a$, in which case $\textbf{x}'$ lies between $\textbf{0}$ and $-\textbf{v}_5$.

Case 13: $\textbf{x}=(-r,s,-a-1,-u)$ where $0 \leq r,s \leq a$ and $0 \leq u \leq a-1$.
Let $\textbf{x}' = \textbf{x} - \textbf{v}_7 = (a-r,s-a, -1,a-1-u)$ lies between $\textbf{0}$ and $-\textbf{v}_3$ unless $s =0$, in which case let $\textbf{x}'' = \textbf{x}' + \textbf{v}_3 = (-1-r, s-1, a-2,-1-u)$ which lies between $\textbf{0}$ and $-\textbf{v}_2$ unless $u=a-1$.
If $s=0$ and $u = a-1$ then let $\textbf{x}''' = \textbf{x}'' + \textbf{v}_2 = (a-r, a+s, 0,a-2-u)$ which lies between $\textbf{0}$ and $\textbf{v}_6$.

Case 14: $\textbf{x}=(-r,a+1,-t,-u)$ where $0 \leq r,t \leq a$ and $0 \leq u \leq a-1$.
Let $\textbf{x}' = \textbf{x} - \textbf{v}_7 = (a-r,1, a-t,a-1-u)$ which lies between $\textbf{0}$ and $-\textbf{v}_5$ unless $u =0$, in which case let $\textbf{x}'' = \textbf{x}' + \textbf{v}_5 = (-r, -a+1, -t-1,1-u)$ which lies between $\textbf{0}$ and $\textbf{v}_8$ unless $r=a$ or $t=a$.
If $r=a$ and $u = 0$ then $\textbf{x}'$ lies between $\textbf{0}$ and $\textbf{v}_4$.
If $0 \leq r \leq a-1$, $t=a$ and $u = 0$ then let $\textbf{x}''' = \textbf{x}'' - \textbf{v}_8 = (a-1-r, 0, a-1-t,-a-u)$ which lies between $\textbf{0}$ and $\textbf{v}_6$.

Case 15: $\textbf{x}=(-a-1,s,-t,-u)$ where $0 \leq s,t \leq a$ and $0 \leq u \leq a-1$.
Let $\textbf{x}' = \textbf{x} - \textbf{v}_7 = (-1,s-a, a-t,a-1-u)$ which lies between $\textbf{0}$ and $-\textbf{v}_6$ unless $t =0$, in which case let $\textbf{x}'' = \textbf{x}' + \textbf{v}_6 = (a-1, s, 1-t,-1-u)$ which lies between $\textbf{0}$ and $-\textbf{v}_8$ unless $s=a$.
If $s=a$ and $t = 0$ then let $\textbf{x}''' = \textbf{x}'' + \textbf{v}_8 = (0, s-a+1, 1-a-t,a-u)$ which lies between $\textbf{0}$ and $\textbf{v}_1$.

This completes the cases for the orthant of $\textbf{v}_7$.

\section*{Orthant of $\textbf{v}_8$, $k \equiv 1$ (mod 2) }
Finally suppose $\textbf{x}$ lies in the orthant of $\textbf{v}_8$ but not between $\textbf{0}$ and $\textbf{v}_8$. Then the first coordinate of $\textbf{x}$ is equal to $-a$ or $-a-1$, or the second is equal to $-a$ or $-a-1$, or the third is equal to $-a-1$. We distinguish seven cases.

Case 1: $\textbf{x}=(-r,-s,-a-1,u)$ where $a \leq r,s \leq a+1$ and $0 \leq u \leq a+1$.
Let $\textbf{x}' = \textbf{x} - \textbf{v}_8 = (a-1-r, a-1-s, -1,u-a-1)$, which lies between $\textbf{0}$ and $\textbf{v}_5$ unless $u \leq 2$, in which case let $\textbf{x}''= \textbf{x}' - \textbf{v}_5 = (2a-1-r,2a-1-s,a,u-3) $ which lies between $\textbf{0}$ and $-\textbf{v}_8$.

Case 2: $\textbf{x}=(-r,-s,-t,u)$ where $a \leq r,s \leq a+1$, $0 \leq t \leq a$ and $0 \leq u \leq a+1$.
Let $\textbf{x}' = \textbf{x} - \textbf{v}_8 = (a-1-r,a-1-s,a-t,u-a-1)$, which lies between $\textbf{0}$ and $-\textbf{v}_2$ unless $t \leq 1$ or $u \leq 1$, in which case let $\textbf{x}''= \textbf{x}' + \textbf{v}_2 = (2a-r,2a-s,2-t,u-2)$. 
If $t \leq 1$ and $u \leq 1$ then $\textbf{x}''$ lies between $\textbf{0}$ and $-\textbf{v}_8$ unless $r=a$ or $s=a$, in which case let $\textbf{x}'''= \textbf{x}'' + \textbf{v}_8 = (a+1-r,a+1-s,2-a-t,u+a-1)$. 
If $t \leq 1$, $u \leq 1$, $r=a$ and $s = a$ then $\textbf{x}'''$ lies between $\textbf{0}$ and $\textbf{v}_2$ unless $t=1$ or $u=1$.
If $r=a$, $s=a$, $t=1$ and $u \leq 1$ then let $\textbf{x}'''' = \textbf{x} - \textbf{v}_5 = (a-r,a-s,a+1-t,a-2+u)$ which lies between $\textbf{0}$ and $\textbf{v}_4$.
If $r=a$, $s=a$, $t=0$ and $u = 1$ then $\textbf{x}$ lies between $\textbf{0}$ and $-\textbf{v}_6$.
If $t \leq 1$, $u \leq 1$, $r=a$ and $s = a+1$ then $\textbf{x}'''$ lies between $\textbf{0}$ and $-\textbf{v}_3$.
If $t \leq 1$, $u \leq 1$, $r=a+1$ and $s = a$ then $\textbf{x}'''$ lies between $\textbf{0}$ and $\textbf{v}_1$.
If $t \leq 1$ and $2 \leq u \leq a+1$ then $\textbf{x}''$ lies between $\textbf{0}$ and $-\textbf{v}_5$ unless $u=a+1$, in which case let $\textbf{x}''''= \textbf{x}'' + \textbf{v}_5 = (a-r,a-s,1-a-t,u-a)$ which lies between  $\textbf{0}$ and $\textbf{v}_8$. 
If $2 \leq t \leq a+1$, and $u \leq 1$ then $\textbf{x}''$ lies between $\textbf{0}$ and $\textbf{v}_6$.

Case 3: $\textbf{x}=(-r, -s,-a-1,u)$ where $a \leq r \leq a+1$, $0 \leq s \leq a-1$ and $0 \leq u \leq a+1$.
Let $\textbf{x}' = \textbf{x} - \textbf{v}_8 = (a-1-r, a-1-s, -1,u-a-1)$, which lies between $\textbf{0}$ and $\textbf{v}_7$ unless $u \leq 1$, in which case let $\textbf{x}''= \textbf{x}' - \textbf{v}_7 = (2a-1-r, -1-s, a-1, u-2)$ which lies between $\textbf{0}$ and $-\textbf{v}_1$.

Case 4: $\textbf{x}=(-r, -s,-a-1,u)$ where $0 \leq r \leq a-1$, $a \leq s \leq a+1$ and $0 \leq u \leq a+1$.
Let $\textbf{x}' = \textbf{x} - \textbf{v}_8 = (a-1-r, a-1-s, -1,u-a-1)$, which lies between $\textbf{0}$ and $-\textbf{v}_4$ unless $u \leq 1$, in which case let $\textbf{x}''= \textbf{x}' + \textbf{v}_4 = (-1-r,2a-1-s,a-1,u-2)$ which lies between $\textbf{0}$ and $\textbf{v}_3$.

Case 5: $\textbf{x}=(-r,-s,-a-1,u)$ where $0 \leq r,s \leq a-1$ and $0 \leq u \leq a+1$.
Let $\textbf{x}' = \textbf{x} - \textbf{v}_8 = (a-1-r, a-1-s, -1,u-a-1)$, which lies between $\textbf{0}$ and $\textbf{v}_6$ unless $u =0$, in which case $\textbf{x}$ lies between $\textbf{0}$ and $\textbf{v}_5$.

Case 6: $\textbf{x}=(-r,-s,-t,u)$ where $0 \leq r \leq a-1$, $a \leq s \leq a+1$, $0 \leq t \leq a$ and $0 \leq u \leq a+1$.
Let $\textbf{x}' = \textbf{x} - \textbf{v}_8 = (a-1-r, a-1-s, a-t, u-a-1)$, which lies between $\textbf{0}$ and $-\textbf{v}_1$ unless $t = 0$ or $u = 0$, in which case let $\textbf{x}''= \textbf{x}' + \textbf{v}_1 = (-r, 2a-2-s, 1-t, u-1)$. 
If $t = 0$ and $u = 0$ then $\textbf{x}$ lies between $\textbf{0}$ and $-\textbf{v}_2$.
If $t =0$ and $1 \leq u \leq a+1$ then $\textbf{x}''$ lies between $\textbf{0}$ and $\textbf{v}_4$ unless $u=a+1$, in which case let $\textbf{x}'''= \textbf{x}'' - \textbf{v}_4 = (a-r, a-2-s, 1-a-t, u-a)$ which lies between  $\textbf{0}$ and $-\textbf{v}_3$.  .
If $1 \leq t \leq a$ and $u =0$ then $\textbf{x}''$ lies between $\textbf{0}$ and $\textbf{v}_7$.

Case 7: $\textbf{x}=(-r,-s,-t,u)$ where $a \leq r \leq a+1$, $0 \leq s \leq a-1$, $0 \leq t \leq a$ and $0 \leq u \leq a+1$.
Let $\textbf{x}' = \textbf{x} - \textbf{v}_8 = (a-1-r, a-1-s, a-t, u-a-1)$, which lies between $\textbf{0}$ and $\textbf{v}_3$ unless $t = 0$ or $u = 0$, in which case let $\textbf{x}''= \textbf{x}' - \textbf{v}_3 = (2a-r, -s, 1-t, u-1)$. 
If $t = 0$ and $u = 0$ then $\textbf{x}$ lies between $\textbf{0}$ and $-\textbf{v}_2$.
If $t =0$ and $1 \leq u \leq a+1$ then $\textbf{x}''$ lies between $\textbf{0}$ and $-\textbf{v}_7$ unless $u=a+1$, in which case let $\textbf{x}'''= \textbf{x}'' + \textbf{v}_7 = (a-r, a-s, 1-a-t, u-a)$ which lies between  $\textbf{0}$ and $\textbf{v}_1$.  .
If $1 \leq t \leq a$ and $u =0$ then $\textbf{x}''$ lies between $\textbf{0}$ and $-\textbf{v}_4$.

This completes the cases for the orthant of $\textbf{v}_8$.

This also completes the proof of the theorem for any $k \equiv 1 \pmod 2$, and therefore for all $k \geq 2$. 
\end{proof}



\end{document}